\RequirePackage{fix-cm}

\documentclass[smallextended]{svjour3}       
\smartqed  

\usepackage{graphicx}
\usepackage{setspace}

\usepackage{float}
\usepackage{breakcites}
\usepackage{epsfig}
\usepackage{amssymb}
\usepackage{amsmath}
\usepackage{hyperref}
\usepackage{verbatim}
\usepackage{subcaption}
\usepackage{tikz}
\usepackage{color}

\usetikzlibrary{arrows}
\usepackage[bottom=0.95in,top=0.95 in]{geometry}
\newcommand\given[1][]{\:#1\vert\:}
\DeclareMathOperator{\E}{\mathbb{E}}
\tikzset{cross/.style={cross out, draw=black, minimum size=2*(#1-\pgflinewidth), inner sep=0pt, outer sep=0pt},
cross/.default={5pt}}
\usetikzlibrary{shapes.misc}
\newcommand{\be}{\begin{eqnarray}}
\newcommand{\ee}{\end{eqnarray}}
\newcommand{\bea}{\begin{eqnarray*}}
\newcommand{\eea}{\end{eqnarray*}}

\newcommand{\s}{\mathbb{S}}
\usepackage{nicefrac}
\numberwithin{equation}{section}

\usepackage{pst-node,pst-plot}
\pstVerb{realtime srand}
\psset{plotpoints=50}

\begin{document}

\title{Approximations for the  boundary crossing probabilities of moving sums of normal random variables}

\titlerunning{Approximations of boundary crossing probabilities}

\author{Jack Noonan    \and Anatoly Zhigljavsky 
}

\institute{Jack Noonan  \at
              School of Mathematics, Cardiff University, Cardiff, CF24 4AG, UK \\
              \email{Noonanj1@cf.ac.uk}           
           \and
           Anatoly Zhigljavsky \at
              School of Mathematics, Cardiff University, Cardiff, CF24 4AG, UK\\
              \email{ZhigljavskyAA@cardiff.ac.uk}
}

\date{Received: date / Accepted: date}

\maketitle

\begin{abstract}
In this paper we study approximations for boundary crossing probabilities for the moving sums of i.i.d. normal random variables. We propose approximating a discrete time problem with a continuous time problem allowing us to apply developed theory for stationary Gaussian processes and to consider a number of approximations (some well known and some not). We bring particular attention to the strong performance of a newly developed approximation that corrects the use of continuous time results in a discrete time setting. Results of extensive numerical comparisons are reported. These results show that the developed approximation is  very accurate even for small window length.


\end{abstract}

\section{Introduction: Statement of the problem}
\label{sec:prob-state}

Let $\varepsilon_1,\varepsilon_2,\ldots$ be a sequence of i.i.d. normal random variables (r.v.) with mean $\mu$ and variance $\sigma^2>0$.
For a fixed positive integer~$L$, the moving sums are defined by
\begin{eqnarray}
S_{n,L}:= \sum_{j=n+1}^{n+L} \varepsilon_j\, \;\; (n=0,1, \ldots) \label{eq:sumsq2}.
\end{eqnarray}
The sequence of the moving sums \eqref{eq:sumsq2} will be denoted by $\s$ so that $\;\s=\{ S_{0,L}, S_{1,L}, \ldots \}$.

The main aim of this paper is development of accurate  approximations
for the following  related characteristics of $\s$
{
(note that for the sake of simplicity of notation we are not indicating the dependence
of these  characteristics on $L$).
}

\begin{enumerate}
  \item[(a)]
 The boundary crossing probability (BCP) for
the maximum of the moving sums:
\be
{\cal P}_{\s}(M,H) := {\rm Pr}\left(\max_{n=0,1,\ldots,M} S_{n,L} \geq H\right),
 \label{eq:prob-S}
\ee
where $M$ is a given positive integer and $H$ is a fixed  threshold. Note that the total number of r.v. $\varepsilon_i$ used in~\eqref{eq:prob-S} is $M+L$ and ${\cal P}_{\s}(M,H) \to 1$ as $M \to \infty$, for all $H$ and $L$.
  \item[(b)] The probability distribution of  the moment of time
  \mbox{$
  \tau_{H}(\s)\!:=\!\min \{ n\!\geq \! 0\!: S_{n,L} \!\geq \! H \}
  $}
   when the sequence $S_{n,L}$  reaches  the threshold $H$ for the first time.
  The BCP ${\cal P}_{\s}(M,H) $, considered as a function of $M$, is the c.d.f. of this probability distribution: ${\cal P}_{\s}(M,H) =
  {\rm Pr}\left(  \tau_{H}(\s) \leq M \right)$.
  \item[(c)]
The average run length (ARL) until $\s$ reaches $H$ for the first time:
\be
\label{arl}
{\rm ARL}_H(\s):= \sum_{n=0}^\infty n {\rm Pr} \{ \tau_H=n\}= \int_0^{\infty}{M d {\cal P}_{\s}(M,H) }  \, .
\ee
\end{enumerate}

Developing accurate approximations for the BCP ${\cal P}_{\s}(M,H) $ and the associated ARL \eqref{arl} for generic parameters $H$, $M$, $L$  is very important in various areas of statistics,  predominantly in applications related to change-point detection; see, for example, papers  \cite{Chu,Glaz2012,MZ2003} and  \cite{glaz2009scan2}.
     The BCP ${\cal P}_{\s}(M,H) $ is an ($M+1$)-dimensional integral and therefore direct evaluation of this BCP is hardly possible even with modern software.

To derive approximations for the  BCP \eqref{eq:prob-S} one can use some generic approximations such as Durbin and Poisson Clumping Heuristic considered below. These approximations, however, are not  accurate especially for small window length $L$; this is demonstrated  below in this paper.
There is, therefore, a need for derivation of specific approximations for the BCP \eqref{eq:prob-S} and the ARL \eqref{arl}. Such a need was well
understood in the statistical community and indeed very accurate approximations for the BCP  and the ARL have been developed in a series of papers
by J. Glaz and coauthors, see for example \cite{Glaz_old} and \cite{Glaz2012}.
 We will call these approximations 'Glaz approximations' by the name of the main author of these papers; they will be formally written down in Sections~\ref{sec:Glaz} and \ref{ARL_section}.

 The accuracy of the approximations developed in the present paper is very similar to the Glaz approximations; this is discussed in Section~\ref{ARL_section}. The methodologies of derivation of the approximations are very different, however.
 The  advantage of the approximation developed in this paper over the Glaz approximation is the fact that our approximation is explicit and hence can be computed instantly; on the other hand,
to compute  the Glaz approximation one needs to numerically approximate $L+1$
and  $2L+1$   dimensional integrals, which is not an easy task even taking into account the fact of  existence of a sophisticated software.

To derive the approximations, in Sections~\ref{sec:Introd} and \ref{sec:da} we have used the methodology developed in \cite[Ch.2,\S 2]{ZhK1988} for continuous-time case, which has to be modified for discrete time. To do this, in  Sections~\ref{discrete_corr_sec}, \ref{Correct_t_gt_1}  and  \ref{4.3} we have revised    and specialized  the approach developed by D.Siegmund in
\cite{Sieg_paper} and other papers.

The paper is structured as follows. In Sections~\ref{Reform}, \ref{C_Diff_section} and \ref{Corr_deriv_T>1} we reformulate the problem and discuss how to approximate our discrete-time problem with a continuous-time problem. Here we state a number of known approximations and derive a new approximation that corrects the use of continuous time results in a discrete time setting; this will be  referred to as the `Corrected Diffusion Approximation' or simply CDA.
In Sections~\ref{sec:sim-study} and \ref{sec:sim-study2} we present results of  large simulation studies evaluating the performance of the considered approximations.
In Section~\ref{ARL_section}, we develop  the  CDA for ${\rm ARL}_H(\s)$ and assess its accuracy.

\section{Boundary crossing probabilities and related characteristics: discrete and continuous time}\label{Reform}
\subsection{Standardisation of the moving sums}

For convenience, we standardise the moving sums $S_{n,L}$ defined in \eqref{eq:sumsq2}.

The first two moments of $S_{n,L}$ are
\be
\label{eq:mean_vars}
\E S_{n,L}= \mu L, \;\;
{\rm var}(S_{n,L})=\displaystyle \sigma^2 L.
\ee
Define the standardized r.v.'s:
\be
\label{eq:def-xi}
\xi_n:= \frac {S_{n,L}-  \E  S_{n,L}}
{\sqrt{{\rm var}(S_{n,L})}}=
\frac { S_{n,L}- \mu L}
{ \sigma\sqrt{  L    }}  \, ,\;\;\mbox{$n=0,1,\ldots\, ,$}
\ee
and denote $\mathbb{X}=\{ \xi_0, \xi_1, \ldots, \}$.
All r.v. $\xi_n$ are $N(0,1)$; that is, they have the probability density function and c.d.f.
\be
\label{eq:phi}
\varphi(x):=\frac{1}{\sqrt{2\pi}}e^{-x^2/2}\,,\; \; \Phi(t):= \int_{-\infty}^t \varphi(x) dx\, .
\ee
Unlike the original r.v. $\varepsilon_i$, the r.v.  $\xi_0, \xi_1, \ldots$  are correlated with correlations depending on~$L$, see Section~\ref{Chapter:correl} below.

The BCP ${\cal P}_{\s}(M,H) $ defined by \eqref{eq:prob-S} is  equal to the BCP
\be
 {\cal P}_{\mathbb{X}}(M,h) := \text{Pr}\left(\max_{n=0,1,\ldots,M} \xi_n \ge h \right),
 \label{eq:prob-xi}
\ee
where
\be
\label{H_h}
H= \mu L +
 \sigma h\sqrt{L} \,\;\mbox{ so that }\; h= \frac{H - \mu L}{\sigma \sqrt{L}}\, .
\ee
Similarly, $\tau_{H}(\s)=\tau_{h}(\mathbb{X})$ and ${\rm ARL}_H(\s)={\rm ARL}_h(\mathbb{X})$.

 In what follows, we derive approximations for \eqref{eq:prob-xi} and hence the distribution of $\tau_{h}(\mathbb{X})$ and ${\rm ARL}_h(\mathbb{X})$.
 These approximations will be
  based on approximating the sequence  $\{\xi_i\}_{i} $ by a continuous time
random process and subsequently correcting the obtained approximations  for discreteness.

\subsection{Correlation between $\xi_{n}$ and $\xi_{n+k}$}
\label{Chapter:correl}

In order to derive our approximations, we will need explicit expressions for the correlations Corr($\xi_{n},\xi_{n+k})$.

\begin{lemma}
\label{lem:corr}
Let  $\xi_{n}$ be as  defined in~\eqref{eq:def-xi}.
Then  ${\rm Corr}(\xi_{0} , \xi_{k})={\rm Corr}(\xi_{n},\xi_{n+k})$ and
\be
\label{eq:correlation}
{\rm Corr}(\xi_0,\xi_k)=
\frac{  \E (\xi_0\xi_k)-( \E \xi_0)^2}{{\rm var}(\xi_0)}=
1-{k}/{L}\, .
\ee
for $0\le k \le L$.
If $k > L$ then ${\rm Corr}(\xi_0,\xi_k)=0$.
\end{lemma}

For a proof, see Appendix A.

\subsection{{Continuous-time (diffusion) approximation}}
\label{sec:cont_time}

For the purpose of approximating the BCP $ {\cal P}_{\mathbb{X}}(M,h) $ and the associated  characteristics introduced in  Introduction,
    we replace the discrete-time process $\xi_0, \ldots, \xi_M$ with a continuous process $\zeta(t)$, $t\in [0,T]$, where  $T=M/L$.
    We do this as follows.

Set $\Delta= 1/L$ and define
$
t_n=n \Delta \in [0, { T}]
\; n=0,1,\ldots ,{M}.
$
Define a piece-wise linear  continuous-time process
${\zeta_t^{(L)}},$ $t \in [0,T]:$
\bea
 {\zeta_t^{(L)}}\!=\!\frac1{\Delta} \left[
(t_{n}-t )\xi_{n-1}\! +\! (t-t_{n-1}) \xi_n \right] \;\;\;{\rm for}
 \;\;t \in [t_{n-1},t_n],\; n=1,\dots,{M}.\;
\eea
By construction, the process ${\zeta_t^{(L)}}$ is such that
${\zeta_{t_n}^{(L)}}=\xi_{n} \; {\rm for } \; n=0,\ldots,{M}$.
Also we have that ${\zeta_t^{(L)}}$ is  a second-order stationary process in the sense that
$ \E \zeta_t^{(L)},\,$ ${\rm var}(\zeta_t^{(L)})$
and the autocorrelation  function
$R_\zeta^{(L)}(t,t+k\Delta)={\rm Corr}( \zeta_t^{(L)}, \zeta_{t+k\Delta}^{(L)})$
do not depend on $t$.

\begin{lemma}
\label{Durbin} Assume $L \to \infty$.
The limiting process $\zeta_t$ = $\lim_{L \rightarrow \infty} \zeta^{(L)}_t$, where $t \in [0,T]$,
is a Gaussian second-order
stationary process with marginal distribution $\zeta_t \sim N(0,1)$ for all $t \in [0,T]$ and autocorrelation function $R_\zeta(t,t+s)=R(s) = \max\{0,\;1\!-\!|s|\}\, $.

\end{lemma}
This lemma is a simple consequence of Lemma~\ref{lem:corr}.

\subsection{Diffusion approximations for the main characteristics of interest}

The above approximation of a discrete-time process $\s$ with a continuous process $\zeta_t,\, t\in [0,T]$, allows us to approximate the  characteristics introduced in  Introduction  by the  continuous-time analogues as follows.

\begin{enumerate}
  \item[(a)]
BCP $
{\cal P}_{\mathbb{X}}(M,h)
$
is approximated by
${\cal P}_{\zeta}(T,h)$, which is the probability of reaching the threshold
$h$ by the process $\zeta_t$ on the interval $[0,{ T}]$:
\be
{\cal P}_{\zeta}(T,h) \!
:=\! {\rm Pr}\left \{\max_{0\leq t\leq { T}}\zeta_t \geq
h\right\}\!  = \!{\rm Pr}\Big\{\zeta_t \! \geq \!
h\;{\rm for\;at\;least\;one\;} t\!\in \! [0,\!{ T}]\Big\} .\;\;
\label{eq:first_pass_prob}
\ee
%
Note that ${\cal P}_{\zeta}(0,h)=1-\Phi(h)>0$.
  \item[(b)] The time moment  $\tau_{H}(\s)=\tau_{h}(\mathbb{X})$ is  approximated
  by  $\tau_{h}(\zeta_t):=\min \{ {t\geq 0}:\;  \zeta_t \geq h \}
  $, which is
   the time moment when the process $\zeta_t$ reaches $h$. The distribution of $\tau_{h}(\zeta_t)$
  has the form:
  \bea
  (1-\Phi(h)) \delta_0 (d s) + q(s,h,\zeta_t) d s \, ,s \geq 0,
  \eea
   where $\delta_0 (d s)$ is the delta-measure concentrated at 0 and
  \begin{equation}
\label{eq:first_pass_time}
q(s,h,\zeta_t)=\frac{d}{ds}{\cal P}_{\zeta}(s,h),\;\;\;0<s<\infty\, .
\end{equation}
The function $q(s,h,\zeta_t)/\Phi(h)$, considered as a function of $s$,  is a probability density function on $(0, \infty)$ since
\bea
\int_0^\infty q(s,h,\zeta_t) ds = 1- {\cal P}_{\zeta}(0,h)= \Phi(h)\, .
\eea

  \item[(c)]
${\rm ARL}_H(\mathbb{X})/L$ is approximated by
\begin{equation}
{\rm ARL}_h(\zeta_t) =  \E (\tau_{h}(\zeta_t))=\int_0^{\infty}{s \,q(s,h,\zeta_t)ds}\, .\label{ARL_form}
\end{equation}

\end{enumerate}
We will call approximations \eqref{eq:first_pass_prob} and  \eqref{ARL_form} diffusion approximations, see Section~\ref{sec:diff_appr}.
Numerical results discussed in Section~\ref{sec:sim-study2} show that if $L$ and $M$ are very large then the diffusion  approximations are  rather accurate.
For not very large  values of $L$ and $M$ these approximations will be much improved with the help of the methodology developed  by D.Siegmund and adapted to our setup in Sections~\ref{discrete_corr_sec} and \ref{Correct_t_gt_1}.

\subsection{Durbin and Poisson Clumping approximations for the BCP ${\cal P}_{\zeta}(T,h)$}

Derivation of the exact formulas for the BCP ${\cal P}_{\zeta}(T,h)$ has been discussed  in several papers including \cite{Mehr,Shepp66,Shepp71,Shepp76,ZhK1988}; exact formulas will be provided in Sections~\ref{Diffusion_section} and \ref{sec:Shepp}.

In this section,  we provide explicit formulas for two simple approximations for the BCP ${\cal P}_{\zeta}(T,h)
$ based on general principles; see also Section~\ref{sec:Glaz} for an approximation specialized for the setup of moving sums.
 We will assess the accuracy of these approximations in Section~\ref{sec:sim-study} and will find that the accuracy of both of them is quite poor. The purpose of including these two approximations into our collection is only to demonstrate that the original problems stated in  Introduction are not
easy and cannot be handled by general-purpose techniques. More sophisticated  techniques using the specificity of the problem should be used, which is exactly what is done in this paper.
The first generic approximation considered is the Durbin approximation which is constructed on the base of  \cite{Durbin} and is explained in Appendix B.\\

\noindent{\bf Approximation 1.} \textit{Durbin approximation for the BCP \eqref{eq:first_pass_prob}:}
$
 {\cal P}_{\zeta}(T,h) \cong hT \, \varphi(h)\, .
$\\

Let us now state the second approximation for the BCP defined in  \eqref{eq:first_pass_prob}, which is the Poisson Clumping Heuristic (PCH) formulated as Lemma~\ref{PCH_lemma} according to \cite{Aldous} p. 81.
\begin{lemma}\label{PCH_lemma}
Let $X(t)$ be a stationary Gaussian process with mean zero and covariance function satisfying $R(t)=1-|t|$ as $t \rightarrow 0$.  Then for large $h$, $T_h = \text{min}\{t: X(t) \ge h  \}$ is approximately exponential with parameter  $h\varphi(h)$.
\end{lemma}

From Lemma~\ref{PCH_lemma} we obtain:\\

\noindent{\bf Approximation 2.} \textit{PCH approximation for BCP \eqref{eq:first_pass_prob}:}
$
{\cal P}_{\zeta}(T,h) \cong 1-\exp(-h\varphi(h) T).
$\\

As can be seen from Fig.~\ref{Quad1} and Fig.~\ref{Quad2} in Section~\ref{sec:sim-study}, Approximations 1 and 2 are poor approximations for ${\cal P}_{\zeta}(T,h)$ and  $ {\cal P}_{\mathbb{X}}(M,h)$ when $M/L\leq 1$; the case $M/L > 1$ is discussed in Section~\ref{sec:sim-study2}.

\section{Diffusion approximation with and without discrete-time correction, $M \leq L$}\label{C_Diff_section}

In this section, we assume $M \leq L$ and hence $T=M/L \leq 1$.  The more complicated case  $M > L$ will be considered in Section~\ref{Corr_deriv_T>1}.

\subsection{Diffusion approximation, formulation} \label{Diffusion_section}

\label{sec:diff_appr}

Here we collect explicit formulas for the BCP  ${\cal P}_{\zeta}(T,h)$; the proofs are given in Section~\ref{sec:Introd}.

We have:
\be
\label{eq:new4}
{\cal P}_{\zeta}(T,h)=1-\Phi^2(h)+\varphi(h) \big[h\Phi(h)
+ \varphi(h) \big]\,, \;\;\;T=1\, ;
\ee

\vspace{-0.4cm}
\begin{equation}\label{eq:new7}
\begin{array}{rcl}
{\cal P}_{\zeta}(T,h)&=&1-\int_{-\infty}^{h}\Phi\left(
\frac{h({ Z}+1)-x(-{ Z}+1)}{2\sqrt{{ Z}}}
\right)\varphi(x)dx+\\[\bigskipamount]
&+& \frac{2\sqrt{Z}}{Z+1} \varphi(h)\left[{h {\sqrt{Z}}}\,\Phi(h\sqrt{{ Z}})
+ \frac1{ \sqrt{2\pi} } (\sqrt{2\pi} \varphi(h))^{Z}\, \right]\, ,\;\;\; 0<T\leq 1\, ,
\end{array}
\end{equation}
where $Z = {T}/({2-T})$. If $T=1$ then   \eqref{eq:new7} simplifies to \eqref{eq:new4}.
We refer to
 the above stated formulas for ${\cal P}_{\zeta}(T,h)$  as Approximation 3 or `Diffusion approximation'.\\

\noindent{\bf Approximation 3.} 
\textit{The Diffusion approximation for the BCP  ${\cal P}_{\mathbb{X}}(M,h) $ defined in \eqref{eq:prob-xi} in case $M\leq L$: formula \eqref{eq:new7} with $T=M/L$; if $M=L$ then \eqref{eq:new7} reduces to
\eqref{eq:new4}.}\\

In Section~\ref{discrete_corr_sec}, we will derive a discrete-time correction for the Diffusion approximation. In order to do this, we need to correct the steps used for deriving \eqref{eq:new7}. This explains that, despite the  formula \eqref{eq:new7} is known, we need to derive it (in order to correct  certain steps of its derivation).
This is done in the next section which follows \cite{ZhK1988}, p.69.

\subsection{Derivation of \eqref{eq:new7}}
\label{sec:Introd}

\subsubsection{Conditioning on the initial value.}
From Lemma~\ref{Durbin}, $\{ \zeta_t$, $t\in [0,\infty) \}$, is a stationary Gaussian process with  mean $\E \zeta_t  = 0 $ and covariance function
$
	\E \zeta_t \zeta_{t+u} = \text{max}\{0,1-|u|\}.
$
By conditioning on the initial state of the process $\zeta_t$, we define
\begin{eqnarray*}
Q_h(T,x_0) := {\rm Pr}\left\{\max_{t\in [0,T]} \zeta_t >h \given \zeta_0 = x_0 \right\}\, .
\end{eqnarray*}
Since $x_0 \sim N(0,1)$ the BCP ${\cal P}_{\zeta}(T,h)$ is
\begin{eqnarray}
\label{three}
{\cal P}_{\zeta}(T,h) &=& \int_{-\infty}^{h}Q_h(T,x_0) \varphi(x_0) dx_0 +1 - \Phi(h) \, ,
\label{6}
\end{eqnarray}
where  $\varphi(\cdot)$ and $\Phi(\cdot)$ are defined in \eqref{eq:phi}.
 In order to proceed we seek an explicit expression for $Q_h(T,x_0)$. We shall firstly discuss a known BCP formula for the Brownian motion before returning to explicit evaluation of  $Q_h(T,x_0)$.

\subsubsection{{Boundary crossing probabilities for the Brownian Motion.}}

Let $W(t)$ be the standard Brownian Motion process on $[0, \infty)$ with $W(0)=0$ and
$
	\E W(t) W(s) =\min(t, s).
$
For given $a,R>0$ and $b \in \mathbb{R}$, define
\begin{equation}\label{BM_notation}
P_W(R; a,b)  := {\rm Pr}\left\{W(t) > a+bt \text{ for at least one }  t\in[0,R] \right\}\, ,
\end{equation}
which is the probability that the Brownian motion $W(t)$ reaches a sloped boundary $a+bt $ within the time interval $[0,R]$.
Using results of \cite{Sieg_paper}, for any $a, R>0$ and any real $b$ we have
\begin{equation}
\label{3.21}
P_W(R; a,b)  = 1 - \Phi\bigg(\frac{bR+a}{\sqrt{R}}\bigg) + e^{-2ab}\Phi\bigg(\frac{bR-a}{\sqrt{R}}\bigg)\, .
\end{equation}
In particular, for $R=1$ we have
\begin{equation}
\label{R_1}
P_W(1; a,b)  = 1 - \Phi(b+a) + e^{-2ab}\Phi(b-a)\, .
\end{equation}
\subsubsection{{Boundary crossing probabilities for $\zeta_t$}.} \label{bcp_for_zeta}
Let $\{\zeta_0(t),t\in[0,\infty)\}$ be a process obtained by considering only the sample functions of $\{ \zeta_t$, $t\in [0,\infty) \}$ which are equal to $x_0$ at $t=0$. For $0 \le t \le 1$, we obtain from \cite{Mehr}, p.520, that $\zeta_0(t)$ can be expressed in terms of the Brownian motion:
\begin{equation}
\zeta_0(t) = (2-t)W( g(t))+ x_0(1-t)   \label{weiner}
\end{equation}
with $g(t)=t/(2-t)$. It then follows from \eqref{weiner} that for $T \le 1$ and $x_0 <h$ we have
\begin{align*}
Q_h(T,x_0)=&  {\rm Pr}\{\zeta_0(t) \ge h \text{ for at least one } t\in[0,T]\}\\
                  =& {\rm Pr}\bigg \{ W(g(t)) \ge \frac{h-x_0(1-t)}{2-t} \text{ for at least one } t\in[0,T]  \bigg \}.
\end{align*}
Noting that   $t={2g(t)}/({1+g(t)})$ we  obtain
\begin{align}
Q_h(T,x_0)=& {\rm Pr}\bigg \{ W(g(t)) \ge \bigg( \frac{(h-x_0)(1+g(t))}{2} \bigg) + x_0g(t) \text{ for at least one } t\in[0,T]  \bigg\} \nonumber\\
=& {\rm Pr}\bigg\{W(t^\prime) \ge \bigg(\frac{h-x_0}{2}\bigg) + t^\prime\bigg(\frac{h+x_0}{2}\bigg) \text{ for at least one } t^\prime \in \bigg[0,\frac{T}{2-T} \bigg] \bigg \} \nonumber\\
=& P_W ( Z;a , b   ) ,\label{Brownian_correction}
\end{align}
where $Z ={T}/({2-T})$, $b= ({h+x_0})/{2}$ and $a = ({h-x_0})/{2}$. Using (\ref{3.21}), we conclude
\begin{gather*}
Q_h(T,x_0)= 1 - \Phi\bigg(\frac{bZ+a}{\sqrt{Z}}\bigg) + e^{-2ab}\Phi\bigg(\frac{bZ-a}{\sqrt{Z}}\bigg)\, .
\end{gather*}
One can then show that by using this explicit form for $Q_h(T,x_0)$ in the integral \eqref{6}, we obtain \eqref{eq:new4} and \eqref{eq:new7}.

It has now become clear how BCP formula \eqref{3.21} for the Brownian motion can be used to obtain \eqref{eq:new4} and \eqref{eq:new7}. To improve the diffusion approximations for discrete time, we aim at correcting the conditional probability $Q_h(T,x_0)$ for discrete time. Because of the relation shown in \eqref{Brownian_correction}, the approach taken in this paper is to correct \eqref{3.21} for discrete time.

\subsection{{Discrete Time Correction }}\label{discrete_corr_sec}

\subsubsection{Discrete time correction for the BCP of cumulative sums.}

Let $X_1,X_2, \ldots $ be i.i.d. $N(0,1)$ r.v's and set $Y_n = X_1+X_2+ \ldots + X_n$.
Consider the sequence of cumulative sums $\{Y_n\}$ and define the stopping time
$
\tau_{Y,a,b} = \inf \{ n \ge 1: Y_n \ge a + bn \}
$
for $a > 0$ and $b \in \mathbb{R}$. Consider the problem of evaluating
\begin{equation}
\label{51}
{\rm Pr}(\tau_{Y,a,b}\le N) = {\rm Pr}(Y_n \ge a+bn \text{ for at least one } n \in \{1,2,\ldots N \}).
\end{equation}

Exact evaluation of \eqref{51} is difficult even if $N$ is not very large but it was accurately approximated by D.Siegmund see e.g. \cite{Sieg_paper} p. 19. Let $W(t)$ be the standard Brownian Motion process on $[0, \infty)$. For $a>0$ and $b \in \mathbb{R}$, define
$
\tau_{W,a,b} = \inf \{ t: W(t) \ge a + bt \}
$
so that
\begin{equation}\label{D_approx}
{\rm Pr}(\tau_{W,a,b}\le N) = P_W(N, a+bt).
\end{equation}

In \cite{Sieg_paper}, \eqref{D_approx} was used to approximate \eqref{51} after translating the barrier $a+bt$ by a suitable scalar $\rho \ge 0$. Specifically, the following approximation has been constructed: \bea
P(\tau_{Y,a,b}\le N) \cong P_W(N, (a+\rho)+bt)\, ,
\eea
where the constant $\rho$ approximates the expected excess of the process $\{Y_n\}$ over the barrier  $a+bt$.
From \cite{Sieg_book} (p. 225)
\begin{equation}
\label{D_rho}
\rho =  - \pi^{-1}\int_{0}^{\infty}\lambda^{-2}  \log\{2(1-\exp(-\lambda^2/2))/ \lambda^2 \} \, d\lambda\, \cong 0.5826.
\end{equation}

Whence, by denoting $\hat{a} = a+\rho$ and recalling \eqref{3.21}, D.Siegmund's formulas of \cite{Sieg_paper} imply the approximation:
\begin{equation*}
{\rm Pr}(\tau_{Y,a,b} \le N) \cong {\rm Pr}(\tau_{W,\hat{a},b} \le N)=  1 - \Phi\bigg(\frac{bN+\hat{a}}{\sqrt{N}}\bigg) + e^{-2\hat{a}b}\Phi\bigg(\frac{bN-\hat{a}}{\sqrt{N}}\bigg).
\end{equation*}

\subsubsection{{Discretized Brownian motion. }}\label{Discretized}

In this section, we modify D.Siegmund arguments discussed in previous section to the case when the r.v. are indexed by points on the uniform grid in an interval and therefore the sequence of cumulative sums compares with a limiting Brownian motion process which lies within this interval.

Assume that $Z>0$ and $M$ is a positive integer. Define $\epsilon$ = $Z/M$ and let
$
t^\prime_n= n\epsilon \in [0, { Z}],
$
$ n=0,1,\ldots ,{M}.$
Let $X_1,X_2, \ldots $ be i.i.d. $N(0,1)$ r.v's and set
$
W(t^\prime_n) = \sqrt{\epsilon}\sum_{i=1}^{{n}}X_i.
$
For $a>0$ and $b \in \mathbb{R}$, define the stopping time
\begin{equation}\label{Brown_stopping}
\tau_{W,a,b} = \inf \{t^\prime_n: W(t^\prime_n) \ge a + bt^\prime_n \}
\end{equation}
and consider the problem of approximating
\begin{equation}
\label{Problem} \!\!\!\!\!
\!{\rm Pr}(\tau_{W,a,b} \le Z) = {\rm Pr}\bigg(W(t^\prime_n)  \ge a \!+\! bt^\prime_n \text{ for at least one } t^\prime_n \in \{0, \epsilon, \ldots,M\epsilon=Z \}\bigg).
\end{equation}

As $M \rightarrow \infty$, the piecewise linear continuous-time process $W^\epsilon(t)$, $t \in [0,Z]$, defined by:
\begin{equation*}
W^\epsilon (t)\!:=\!\frac1{\epsilon} \left[
(t^\prime_{n}-t )W(t^\prime_{n-1})\! +\! (t-t^\prime_{n-1}) W(t^\prime_{n})\right] \;\;\;{\rm for}
 \;\;t \in [t^\prime_{n-1},t^\prime_n],\; n=1,\dots,{M},\;
\end{equation*}
converges to the Brownian motion on $[0,Z]$. For this reason, we refer to the sequence $\{ W(t^\prime_1),\ldots W(t^\prime_M),\}$  as discretized Brownian motion. We make the following connection between $W(t^\prime_n)$ and the random walk $Y_n$:

\begin{equation*}
W(t^\prime_n)=\sqrt{\epsilon} \,Y_n = {Y_n}/{\sqrt{M/Z}}\,\,, n=1,2,\ldots M.
\end{equation*}
Then by using \eqref{D_rho}, we approximate the expected excess over the boundary for the process $W(t^\prime_n)$  by
$
\rho_{M/Z}  =  {0.5826}/{\sqrt{M/Z}}\, .
$

Thus, using the same methodology as D.Siegmund, in order to obtain an accurate approximation for \eqref{Problem}, we translate the barrier $a+bt$ by the discrete time correction factor $\rho_{M/Z}$ and apply \eqref{3.21}. By denoting $ \hat{a} = a+\rho_{M/Z}$, we obtain the approximation to \eqref{Problem}:
\begin{equation}\label{brown_approx}
{\rm Pr}(\tau_{W,a,b} \le Z) \cong  1 - \Phi\bigg(\frac{bZ+\hat{a}}{\sqrt{Z}}\bigg) + e^{-2\hat{a}b}\Phi\bigg(\frac{bZ-\hat{a}}{\sqrt{Z}}\bigg).
\end{equation}

\subsubsection{{Corrected Diffusion Approximation.}}

 Let  {$Q_{h,\rho}(M,x_0)$} denote the discrete time corrected equivalent of $Q_{h}(T,x_0)$, where $T=M/L\leq 1$. Using \eqref{brown_approx}
  and the relation shown in \eqref{Brownian_correction},

\begin{align}\label{discrete_conditional}
Q_{h,\rho}(M,x_0) =1 - \Phi\bigg(\frac{bZ+\hat{a}}{\sqrt{Z}}\bigg) + e^{-2\hat{a}b}\Phi\bigg(\frac{bZ-\hat{a}}{\sqrt{Z}}\bigg)
\end{align}
with
\begin{equation*}
 T=\frac{M}L, \,\,\, Z = \frac{T}{2-T}, \,\,\, \hat{a}=\frac{h-x_0}{2} + \rho_{M/Z}, \,\,\, b = \frac{h+x_0}{2}, \;\;\rho_{M/Z}= \frac{0.5826}{\sqrt{M/Z}}.
 \end{equation*}

Using $Q_{h,\rho}(M,x_0)$ in \eqref{three}, the equivalent probability ${\cal P}_{\zeta}(T,h)$ after correction for discrete time will be denoted by ${\cal P}_{\zeta,\rho}(M,h)$.

\noindent{\bf Approximation 4.} \textit{For $M \le L$  (that is, $T \le 1$), the CDA for the BCP \eqref{eq:prob-xi} is given by}

\begin{equation}\label{Corrected_Diffusion}
{\cal P}_{\mathbb{X}}(M,h)  \cong {\cal P}_{\zeta,\rho}(M,{h}) :=  \int_{-\infty}^{h} Q_{h,\rho}(M,x_0) \varphi(x_0) dx_0 +1- \Phi(h)\, ,
\end{equation}
where $Q_{h,\rho}(M,x_0)$ is given in \eqref{discrete_conditional}.

{For  $M=L$ we have $T=Z=1$ and  the CDA
 ${\cal P}_{\zeta,\rho}(M,{h})$ can be explicitly evaluated}:
\begin{equation}\label{Corrected_Diffusion_explicit}
{\cal P}_{\zeta,\rho}(L,{h})= 1 - \Phi(h+\rho_{L})\,\Phi(h) + \frac{\varphi(h+\rho_{L})}{\rho_{L}}\Phi(h) - \frac{\varphi(h)e^{-2h\rho_{L}}}{\rho_{L}}\Phi(h-\rho_{L})\, ,
\end{equation}
where $\rho_{L}=  0.5826/ \sqrt{L} $.
For a proof of \eqref{Corrected_Diffusion_explicit}, see Appendix C.


\subsection{Simulation study, $T \le 1$}
\label{sec:sim-study}

\begin{figure}[h]
\begin{center}
 \includegraphics[width=0.5\textwidth]{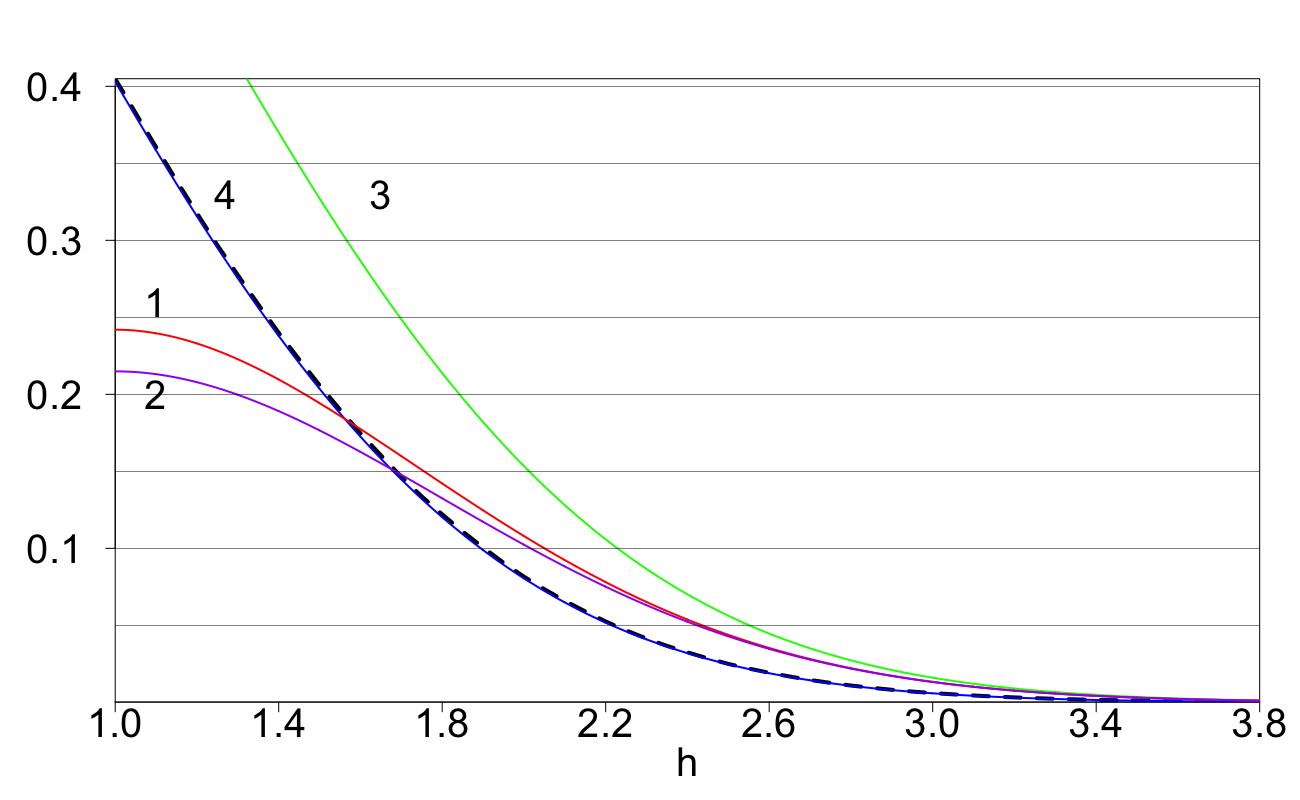}\includegraphics[width=0.5\textwidth]{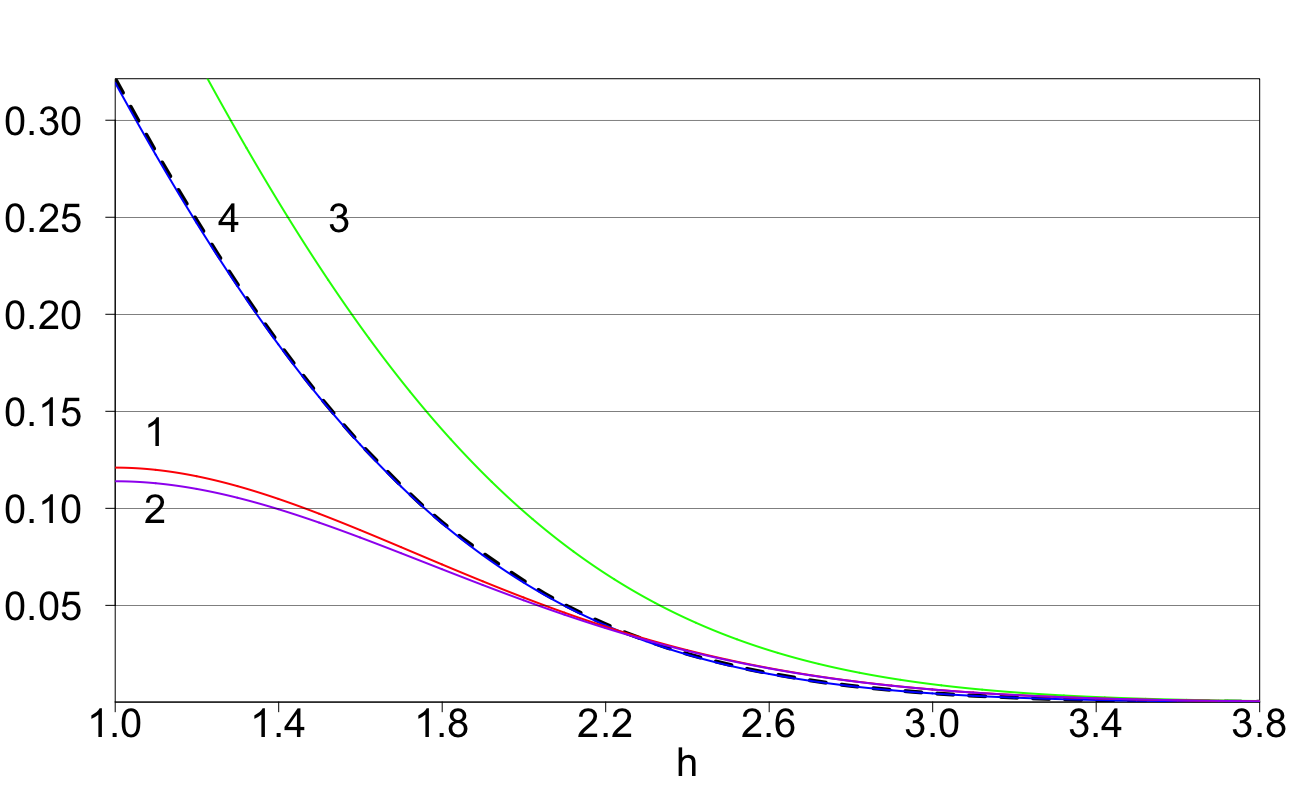}
\end{center}
\caption{Empirical probabilities of reaching the barrier $h$ and four approximations.
Left: $L=5$, $M=5$, $T=1$. Right: $L=10$, $M=5$, $T=1/2$. }
\label{Quad1}
\end{figure}

\begin{figure}[h]
\begin{center}
 \includegraphics[width=0.5\textwidth]{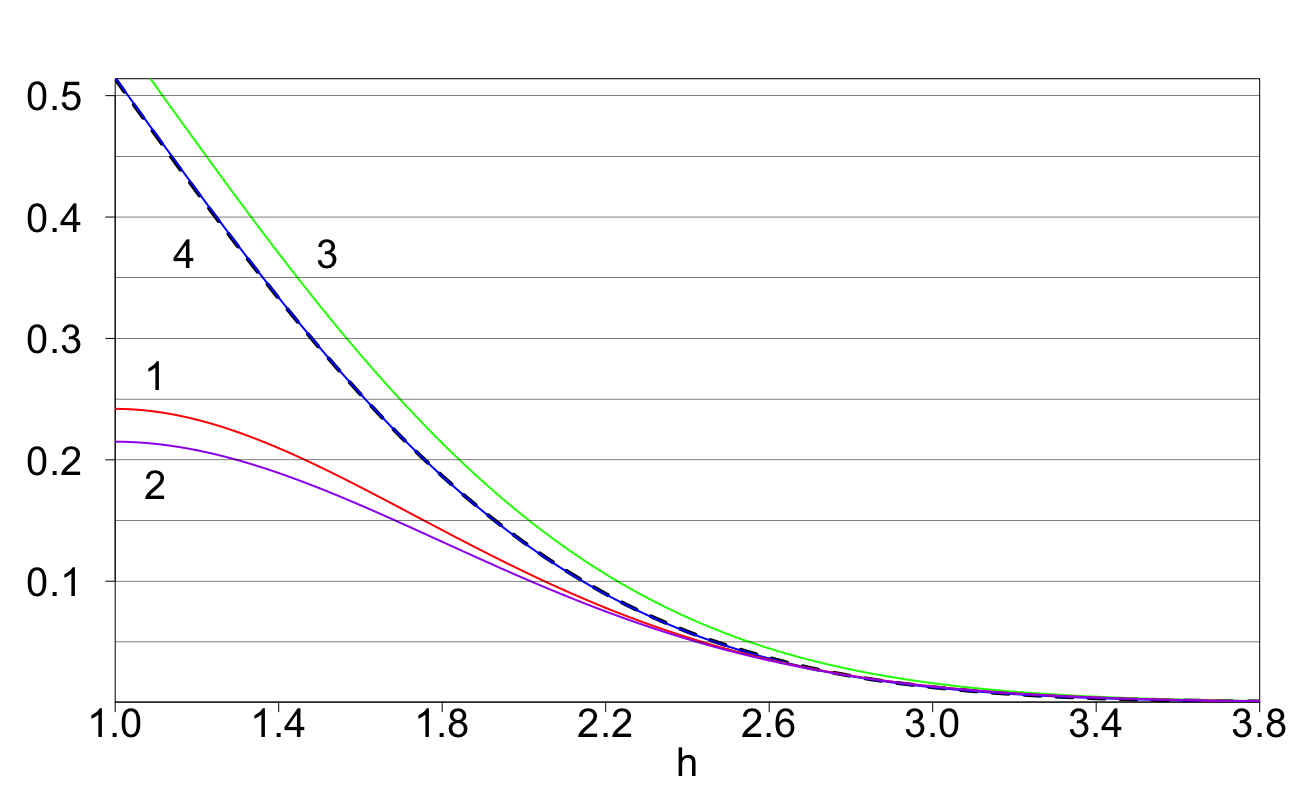}\includegraphics[width=0.5\textwidth]{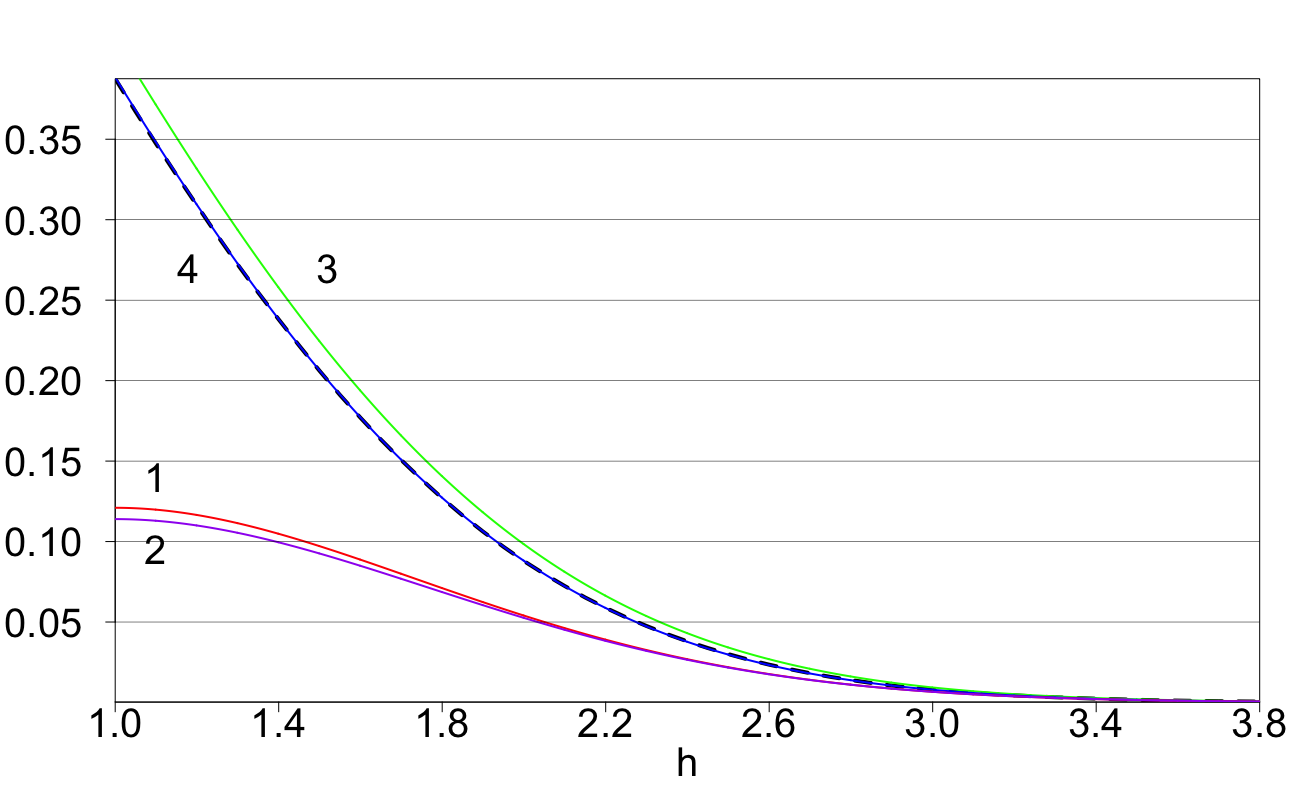}
\end{center}
\caption{Empirical probabilities of reaching the barrier $h$ and four approximations. Left: $L=100$, $M=100$, $T=1$. Right: $L=200$, $M=100$, $T=1/2$. }
\label{Quad2}
\end{figure}

In this section we study the quality of the Durbin (Approximation 1), PCH (Approximation 2), Diffusion (Approximation 3) approximations and the CDA (Approximation 4) for the BCP
${\cal P}_{\mathbb{X}}(M,h)$, defined in  \eqref{eq:prob-xi},  when $M\le L$ (that is, $T\le1$). Without loss of generality, $\varepsilon_j$ in \eqref{eq:sumsq2} are normal r.v.'s with mean $0$ and variance $1$. In Figures~\ref{Quad1}--\ref{Quad2}, the black dashed line corresponds to the empirical values of the BCP ${\cal P}_{\mathbb{X}}(M,h)$ defined by \eqref{eq:prob-xi} computed from 100\,000 simulations  with different values of $L$ and $M$ (for given $L$ and $M$, we simulate $L+M$ normal random variables 100\,000 times). For $j=1,\ldots, 4$, the number~$j$ next to a line corresponds to Approximation $j$. The axis are: the $x$-axis shows the value of the normalized barrier $h$, see \eqref{H_h}; the $y$-axis denotes the probabilities of reaching the barrier. The graphs, therefore,  show the empirical probabilities of reaching the barrier $h$ (for the
dashed line) and values of considered approximations for these probabilities.

In Table \ref{Relative_error_table}, we display the relative error of the  CDA with respect to  the empirical BCP ${\cal P}_{\mathbb{X}}(M,h)$ for all considered parameter choices. Numerical study of this section shows that in the case $T\leq 1$,  the accuracy of the CDA (Approximation 4) is excellent, even for rather small $L$ and $M$. At the same time, the Durbin, PCH and Diffusion approximations  are generally poor (note however that the accuracy of the Diffusion approximation improves as $L$ increases). The discrete time correction factor brings a huge improvement to the Diffusion approximation resulting in a very small relative errors shown in Table~\ref{Relative_error_table}.

\vspace{-0.2cm}

{
\begin{table}[h]
\centering
\caption{Relative error of the CDA with respect to the empirical BCP (in percent)}      \label{Relative_error_table}
\begin{tabular}{|c||c|c|c|c|c|c|}
\hline
BCP & $L\!=\!5,M\!\!=\!5$& $L\!=\!10,M\!\!=\!5$  & $L\!=\!100,M\!\!=\!100$& $L\!=\!200,M\!\!=\!100$\\
\hline
0.05 &  0.225 \% &  0.238 \% &   0.041   \%&  0.132   \%  \\
0.10 &  0.316 \%  & 0.284  \% &  0.093    \%&  0.103  \% \\
0.15 &  0.474 \% &  0.326 \% &   0.155    \%&  0.059  \% \\
0.20 &  0.390 \% &  0.296 \% &   0.228   \%&   0.101  \%  \\
\hline
\end{tabular}
\end{table}
}

\section{Approximations for the BCP in continuous and discrete time; $M > L$}\label{Corr_deriv_T>1}

In this section, we assume  $M>L$ and thus $T>1$.

\subsection{Exact formulas for the continuous-time BCP ${\cal P}_{\zeta}(T,h)$}
\label{sec:Shepp}

For $T>1$, the exact formulas for the BCP ${\cal P}_{\zeta}(T,h)$, the continuous-time case, are  complicated.
 If $T>1$ is an integer then the results of  \cite{Shepp71} imply
\begin{equation}\label{shepp_form}
{\cal P}_{\zeta}(T,h) = 1- \int_{-\infty}^{h}\int_{D_x} \det|\varphi(y_i - y_{j+1} + h){|^T} \hspace{-0.31cm}|_{i,j=0} \, dy_2\ldots dy_{T+1}\,dx,
\end{equation}
where $y_0= 0, y_1=h-x,$
$
D_x =  \{y_2, \dots , y_{T+1} \given h-x < y_2 < y_3 < \ldots < y_{T+1}   \}
$.
If $T>1$ and is not an integer, the exact  formula for ${\cal P}_{\zeta}(T,h)$ is even more complex, see \cite{Shepp71}.

We are not using the exact formulas  for
${\cal P}_{\zeta}(T,h)$ in the case $T>1$ in our approximations for the following two reasons: (a)  the formulas are complicated and  (b) we do not know, yet, how to correct \eqref{shepp_form} for discrete time.
Instead, we  shall derive an approximation for computing BCP ${\cal P}_{\zeta}(T,h)$, which we will also call `Diffusion approximation', and then correct it for discrete time.

\subsection{{A Diffusion approximation when $T>1$}}
\label{sec:da}

To proceed, we need  the following result for the standard Brownian motion process $W(t)$.

\begin{lemma} {\rm (\cite{Harr}, Corollary on p.12 \rm )}
Let  $W(t)$ be the standard Brownian motion on $[0, R]$ with $W(0)=0$ and $ \E W(t)W(s)= \min\{t,s\}$. Let $W_\mu(t)= \mu t + W(t)$ be the Brownian motion with drift $\mu t$. Then, for any  $y>0$ and $R>0$,
\bea \label{eq:absorb}\!\!
F_{R,\mu}(z,y)\!:=\!{\rm Pr} \{ W_\mu(R)\leq z, \sup_{t \in [0,R]} W_\mu(t) \!\leq \! y \}\!= \!\Phi\left(\frac{z\!-\!\mu R}{\sqrt{R}}\right)\!-\!  e^{2y\mu}\Phi\left(\frac{z\!-\!\mu R\!-\!2y}{\sqrt{R}}\right). \;
\eea
\end{lemma}
Similarly to formula (2) on p. 11 in \cite{Harr} we can write the above formula  in the form
\bea
{\rm Pr} \{ W_\mu(R) \in dz, \sup_{t \in [0,R]} W_\mu(t) \leq y \} = f_{R,\mu,y}(z)\, dz,
\eea
where
$$
 f_{R,\mu,y}(z) = \frac{\partial F_{R,\mu}(z,y)}{\partial z} \! =\! \frac{1}{\sqrt{2\pi R}}\bigg\{\exp \bigg[-\frac{(z-\mu R)^2}{2R} \bigg] \!-\! \exp\bigg[2y\mu \!-\! \frac{(z\!-\!2y\!-\!\mu R)^2}{2R} \bigg] \bigg\} , \; z<y \,.
$$
From the definition of $F_{R,\mu}(z,y)$,
\be
\label{bound_prob}
\int_{-\infty}^y f_{R,\mu,y} (z)  \, dz= {\rm Pr} \left \{ \sup_{t \in [0,R]} W_\mu(t) \leq y \right \}= 1-P_W(R;a,b)\, ,
\ee
where $a=y$, $b=-\mu$ and $P_W(R;a,b)$ is defined in \eqref{BM_notation}.

We can reformulate the above results stated for $W_\mu(t)$ as results for the standard Brownian motion process  $W(t)$ with no drift. Set
\begin{equation}
\label{2.17}
p_W(x;R,a,b):=f_{R,-b,a}( x\!-\!bR)=
\begin{cases}
\frac{1}{\sqrt{2\pi R}}\left[ e^{-x^2/2R} \! - \! e^{-2ab \! - \! {(x\! -\! 2a)^2/2R}  }    \right], &  x \! \le a\!+\!bR\\
0, &  x>a\!+\!bR.
\end{cases}
\end{equation}

We will call $p_W(x;R,a,b)$ in  (\ref{2.17})
`the non-normalised density
 of $W(R)$ under the condition $W(t) < a+bt$ for all $t \in [0,R]$'.
 In view of (\ref{bound_prob}),
 $\int p_W(x;R,a,b)dx =1- P_W(R;a,b).
$

Let us now show how to apply these results for construction of approximations for the BCP ${\cal P}_{\zeta}(T,h)$, where $\zeta_t$ is the process defined in Lemma~\ref{Durbin}. The direct relation between the process $\zeta_t$, $t \in [0,1]$, and the standard Brownian motion   is given in   \eqref{weiner}.

 Let $V \in (0,1] $. From  \eqref{Brownian_correction},  the  conditional probability that $\zeta_t<h$ for all $t\in(0,V]$ given $\zeta_0=x_0<h$ is
\begin{align*}
{\rm Pr}(\zeta_t < h \text{ for all } t\in [0,V]\given \zeta_0=x_0) = 1-P_W(U;a,b)\, ,
\end{align*}
where $U = {V}/{(2-V)}$, $a={(h-x_0)}/{2}$ and $b = {(h+x_0)}/{2}$. By substituting  these particular $a$ and $b$ into (\ref{2.17}), we obtain
that  the non-normalised density of the r.v. $\zeta_V$ conditioned on $\zeta_0 = x_0<h$ and $\zeta_t < h$ for all $t \in (0,V]$ is
\be
\label{2.19}
p_{h,V}(x\given x_0) = \sqrt{\frac{2-V}{2 \pi V}}\bigg\{ \exp\bigg[-\frac{(2-V)x^2}{2V}\bigg] - \exp\bigg[\frac{(2-V)(x-h+x_0)^2}{2V}  \bigg]   \bigg\}.
\ee
In the most important special case  $V=1$,
the non-normalised density of the r.v. $\zeta_1$ conditioned on $\zeta_0 = x_0$ and $\zeta_t < h$ for all $t \in (0,1]$ is
\begin{equation}
\label{2.20}
p_h(x\given x_0) =
\left\{
  \begin{array}{lc}
    \varphi(x)\left[1 - \exp\{-(h-x)(h-x_0)\} \right], & \text{     for $x<h$} \\
    0 & \text{     for $x\geq h$}, \\
  \end{array}
\right.
 \end{equation}
where $\varphi(\cdot)$ is defined in \eqref{eq:phi} and $\int_{-\infty}^{h}p_h(x\given x_0) \, dx = 1-P_W(1;\frac{h-x_0}{2},\frac{h+x_0}{2})$, where we have used \eqref{R_1} to get the final expression.

 Since $\zeta_0$ is $N(0,1)$, the density  of $\zeta_0$ conditioned on $\zeta_0 <h$ is ${p}_0(x)= \phi(x)/\Phi(h), \;x<h$. Averaging over $\zeta_0<h$, the non-normalized density of the r.v. $\zeta_1$ under the conditions $\zeta_t<h$ for all $t \in[0,1]$ is:
\begin{equation}
\label{2.22}
\tilde{p}_1(x) = \int_{-\infty}^{h}p_h(x \given x_0){p}_0(x_0)dx_0, \text{      for $x<h$}
\end{equation}
with
\be
\label{eq:c1}
c_1=\int_{-\infty}^{h}\tilde{p}_1(x)dx =[1-{\cal P}_{\zeta}(1,h)]/\Phi(h)\, .
\ee

Denote by ${p}_1(x)=\tilde{p}_1(x)/c_1, \; x<h,$  the normalized density of
 $\zeta_1$ under the condition $\zeta_t<h$ for all $t \in[0,1]$.

For any integer $i \geq 1 $, the densities of $\zeta_i$ and $\zeta_{i-1}$ under the condition that  $\zeta_t$ does not reach  $h$ in $(i-1,i]$ and $[0,i-1]$ respectively can be connected in the same way as for the interval $[0,1]$ (note, however, that these are only approximations as the process $\zeta_t$ is not Markovian).
Assume that
${p}_{i-1}(x)$ is the normalized density of
 $\zeta_{i-1}$ under the condition $\zeta_t<h$ for all $t \in[0,i-1]$.
Define
 \begin{equation}
\label{2.24}
\tilde{p}_i(x) = \int_{-\infty}^{h}p_h(x\given y)p_{i-1}(y)dy, \text{    for $x<h$}\,.
\end{equation}
We  call it  the non-normalized density of a r.v. $\zeta_{i}$   under the conditions $\zeta_{i-1} \sim \tilde{p}_{i-1}(x)$ and $\zeta_t<h$ for all $t \in[i-1,i]$. We then define ${p}_i(x)=\tilde{p}_i(x)/c_i, \; x<h,$
where $c_i=\int_{-\infty}^{h}\tilde{p}_i(x)dx$.
%

If $T$ is large, calculation of the densities $p_i(x)$ $(i\leq T)$ in such an iterative manner is cumbersome. 
We then replace formula (\ref{2.24}) with
\begin{equation}
\label{2.26}
\tilde{p}_i(x) = \int_{-\infty}^{h}p_h(x\given y)p(y)dy, \text{    for $x<h$},
\end{equation}
where $p(x)$ is an eigenfunction of the integral operator with kernel (\ref{2.20}) corresponding to the maximum eigenvalue $\lambda$:
\be
\lambda p(x) = \int_{-\infty}^{h}p(y)p_h(x\given y) dy, \text{  } x<h\, . \label{2.27}
\ee
This eigenfunction $p(x)$ is a probability density on $(-\infty,h]$ with $p(x)> 0$ for all $x \in (-\infty,h)$ and
$
\int_{-\infty}^{h}p(x)dx = 1 \, .\label{2.28}
$
Moreover, the maximum eigenvalue $\lambda$ of the operator with kernel $K(x,y)=p_h(x|y)$ is simple and positive.
The fact that such maximum eigenvalue $\lambda$ is  simple and real (and hence positive) and the eigenfunction $p(x)$ can be chosen as  a probability density follows from the Ruelle-Krasnoselskii-Perron-Frobenius theory of bounded linear positive operators, see e.g. Theorem XIII.43 in \cite{ReedSimon}.

Using  (\ref{2.26}) and (\ref{2.27}), we  derive recursively:
\begin{gather*}
{\cal P}_{\zeta}(i,h) \simeq  {\cal P}_{\zeta}(i-1,h) + (1-{\cal P}_{\zeta}(i-1,h))(1-\lambda)\, ; \;\;{i=2,3,\ldots}
\end{gather*}
By induction,  for an integer $T\geq 2$ we then have
\be
\label{appr_main}
 {\cal P}_{\zeta}(T,h) \simeq \,  {\cal P}_{\zeta}(1,h) + (1- {\cal P}_{\zeta}(1,h))\sum_{j=0}^{T-2} \lambda^j(1-\lambda)
                                        = \, 1- (1-{\cal P}_{\zeta}(1,h) ) \lambda^{T-1} \, .
\ee
The approximation \eqref{appr_main} can  be used for non-integer $T$. We can also use a minor adjustment to this approximation using the maximal eigenvalues of the kernel \eqref{2.19} in addition to $\lambda$;  this is much harder but the benefits of this are minuscule.\\

{\bf Approximation 5.} (Diffusion approximation for BCP ${\cal P}_{\zeta}(T,h) $ when $T>1$). {\it Use \eqref{appr_main}, where $\lambda$
is the maximal eigenvalue of the integral operator with kernel $p_h(x\given y)$ defined in \eqref{2.20}.}\\

The BCP ${\cal P}_{\zeta}(1,h) $  can be computed by \eqref{eq:new4}. In the next section we make a correction to Approximation 5 adjusted for discrete time. Approximations for continuous-time $\lambda$, required in Approximation 5, are obtained from formulas of that section when  $\rho \rightarrow 0$.

\subsection{{Corrected Diffusion approximation, $T>1$}}\label{Correct_t_gt_1}
There are two components of the Diffusion approximation (Approximation 5) that can and should be corrected for discrete time. These are:
(a) the BCP ${\cal P}_{\zeta}(1,h)$, and
(b) the kernel of the integral operator defined  in \eqref{2.20} for the continuous-time case (and hence $\lambda$, the maximum eigenvalue).

Discrete-time corrections for  the BCP ${\cal P}_{\zeta}(1,h)$ have been discussed above in Section~\ref{discrete_corr_sec}; we will return to this at the end of Section~\ref{4.3}.

Moving to  (b), recall \eqref{2.27} which states that $\lambda$ is the maximal eigenvalue satisfying \eqref{2.27},
 the corresponding eigenfunction $p(x)$  is a probability density function on $(-\infty, h)$, and $p_h(x\given y)$ is given in \eqref{2.20}. Recall that $p_h(x \given x_0)$ is the density of the random variable $x=\zeta_1$ under the conditions $\zeta_0 = x_0$ and $\zeta_t < h$ for all $t \in [0,1]$. We shall now discuss  how to correct the kernel $p_h(x \given x_0)$ for discrete time.

As shown in Section~\ref{Discretized}, BCP for the discretized Brownian motion process $W(t_n)$ can be approximated using the BCP formula for the Brownian motion incorporating a discrete time correction factor. Recalling the stopping time $\tau_{W,a,b}$ defined in \eqref{Brown_stopping},  for $Z=1$ we obtain from \eqref{brown_approx} the approximation:
\begin{equation*}
{\rm Pr}(\tau_{W,a,b} \le 1) \cong  1 - \Phi(b+\hat{a}) + e^{-2\hat{a}b}\Phi(b-\hat{a})\, ,
\end{equation*}
where $\hat{a} = a+\delta$ and, since $L=M$, we will use
\begin{equation}
\label{delta}
\delta  := \rho_L = {0.5826}/{\sqrt{L}}\, .
\end{equation}
In view of \eqref{brown_approx}, for the  Brownian motion $W(t)$ (with $W(0)=0$) considered on  $[0,1]$, the non-normalised density of $W(1)$ under the condition $W(t_n)<a+bt_n$  for all $t_n = n/L $ ($n=1, \ldots, L$) can be approximated by
\begin{equation*}
\psi_\delta(x)=
\begin{cases}
\frac{1}{\sqrt{2\pi}}\bigg\{ e^\frac{-x^2}{2} - \exp\bigg[-2\hat{a}b - \frac{(x-2\hat{a})^2}{2} \bigg]     \bigg\} & \text{for } x \le a+b\\
0 & \text{for } x>a+b\, ,
\end{cases}
\end{equation*}
where $\hat{a} = a+\delta$ and $\delta$ is given by \eqref{delta}.
Thus, the discrete-time equivalent of $p_h(x \given x_0)$, denoted by $p_{h,\delta}(x \given x_0)$, is:
\begin{equation}\label{Discrete_kernel}
p_{h,\delta}(x\given x_0) = \varphi(x)(1 - \exp[-(h-x)(h-x_0) - \delta(3h-2x -x_0+2\delta)] ), \text{     for $x<h$}.
\end{equation}
If $\rho=\delta=0$, we clearly obtain $p_{h,0}(x\given x_0)=p_{h}(x\given x_0)$, for all $h,x,x_0$.

Denote by $\lambda_{\delta}$  the maximum eigenvalue associated with the kernel ${p}_{h,\delta}(x\given x_0)$ given by \eqref{Discrete_kernel}. This means that  $\lambda_\delta$ satisfies
\be
\label{eq:lambdaD}
\lambda_\delta p(x) = \int_{-\infty}^{h}p(y)p_{h,\delta}(x\given y) dy, \text{  } x<h
\ee
for some  eigenfunction $p(x)$ which can be assumed to be a probability density function on $(-\infty,h)$. Similar to $\lambda$ in \eqref{2.27}, $\lambda_\delta$ is positive and uniquely defined.

\subsection{Approximating $\lambda_\delta$ and $p(x)$ in \eqref{eq:lambdaD}}\label{4.3}
\begin{figure}[!h]
\begin{center}
 \begin{subfigure}{0.48\linewidth} \centering
     \includegraphics[scale=0.22]{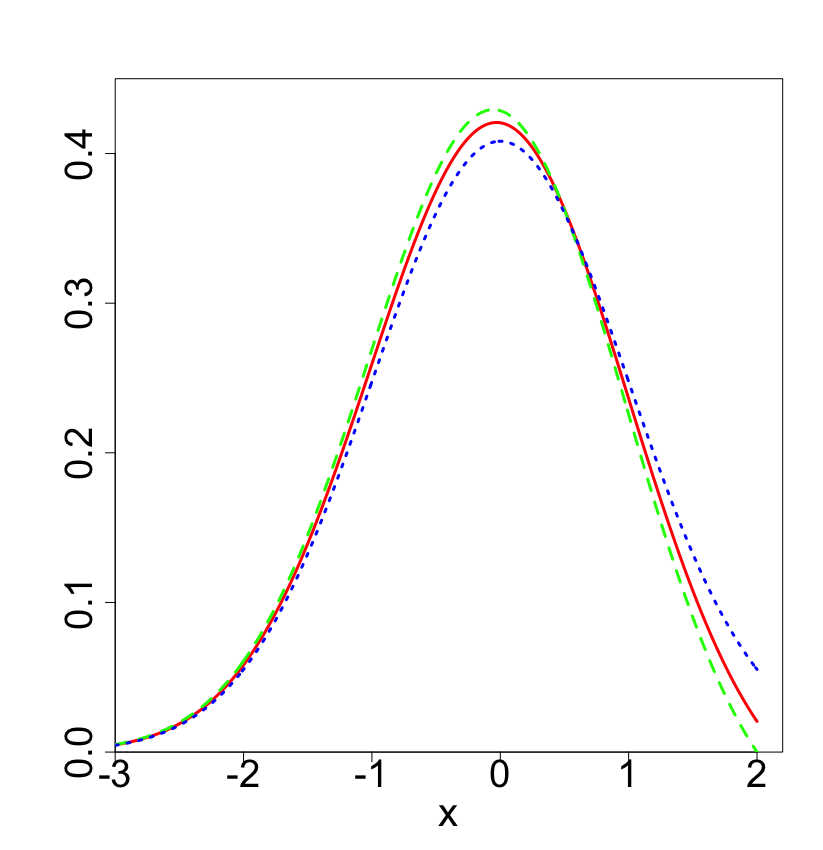}
     \caption{Densities $p_0(x)$ (dotted blue), $p_1(x)$ (solid red) and the uncorrected $p_1(x)$ with $\delta=0$ (dashed green); $L=10$ and $h=2$}\label{fig:densities}
   \end{subfigure} \;
   \begin{subfigure}{0.48\linewidth} \centering
     \includegraphics[scale=0.22]{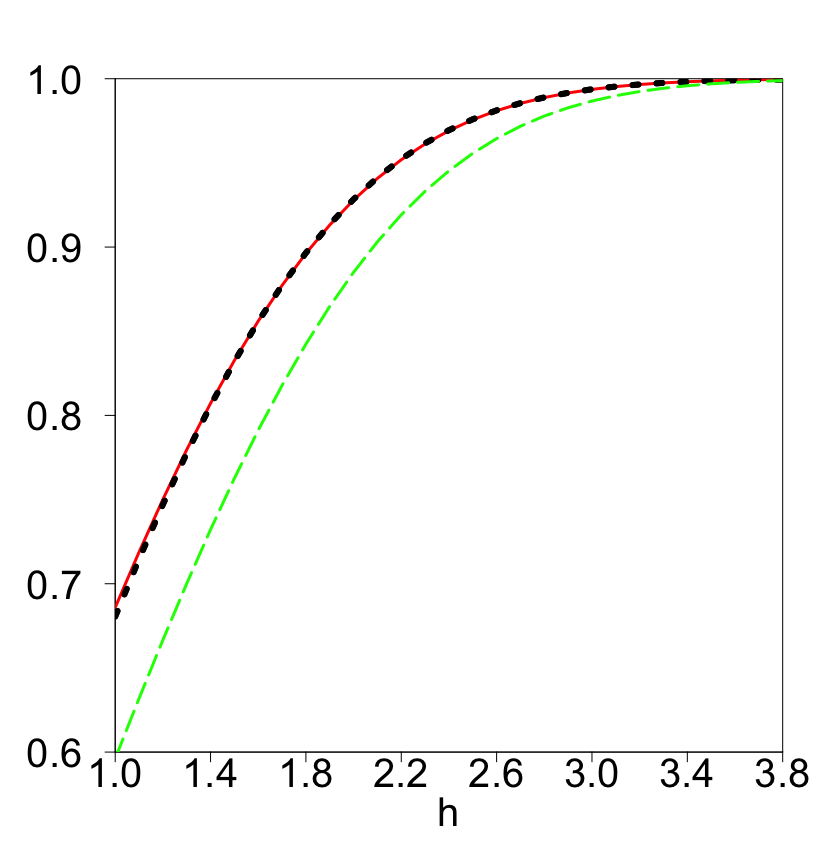}
     \caption{Values of $\hat{\lambda}_\delta$  (solid red),
        ${\lambda}_\delta$ (dotted black) and ${\lambda}_0$ (dashed green) obtained using the Gauss-Legendre quadrature; $L=10$ and different $h$.} \label{fig:lambda}
   \end{subfigure}
\end{center}
\caption{ }
\label{d_and_l}
\end{figure}

Similarly to \eqref{2.22},  \eqref{eq:c1} and
\eqref{2.24} we define $p_0(x)=\phi(x)/\Phi(h)$,
$\tilde{p}_i(x) = \int_{-\infty}^{h}p_{h,\delta}(x\given y)p_{i-1}(y)dy$   ($x<h$),
 $c_i=\int_{-\infty}^{h}\tilde{p}_i(x)dx$ and ${p}_i(x) = \tilde{p}_i(x)/c_i$ for $i=1,2$. From \eqref{eq:lambdaD}, $\lambda_\delta$  = $\lim_{i\to\infty} c_i$.
 By performing integration we obtain
\bea
\tilde{p}_1(x) = \varphi(x) - \frac{\varphi(h)}{\Phi(h)}e^{-2\delta h   - 3\delta^2 /2 + \delta x }\Phi(x-\delta)\, , \;\; x < h\, ,
\eea
\bea
c_1=
\Phi \left( h \right) -{\frac {\kappa_{h,\delta}  }{\Phi \left( h \right) }}\;\;\;{\rm with}\;\; \kappa_{h,\delta}=\frac1{\delta} \varphi(h) \left[ {{\rm e}^{-\delta h-3{\delta}^{2}/2}}
\Phi \left( h-\delta \right) -{{\rm e}^{-2\delta h}}\Phi \left( h-2\,
\delta \right) \right]\, ,
\eea
\bea
\tilde{p}_2(x)= \varphi(x)+ \frac{\varphi(h)}{c_1}\bigg[   \frac{\Phi(h-\delta)e^{- 3\delta h -7\delta^2/2  +2\delta x}\varphi(x) - \varphi(h)e^{\delta^2/2-2\delta h - \delta x}\Phi(x-3\delta) }{\Phi(h)(h+2\delta -x)}
 \\
-\, e^{-2\delta h -3\delta^2/2+\delta x}\Phi(x-\delta) \bigg] , \;\; x < h\, .
\eea
We were unable to compute  $c_2$ and the densities $p_i(x)$ with $i\geq 2$ analytically. However, numerical computations show that the density $p_1(x)$ is visually indistinguishable from $p_i(x)$ for $i>1$ and hence from $p(x)$, the solution of \eqref{eq:lambdaD}. Thus we approximate $\lambda_\delta$ in \eqref{eq:lambdaD} by
\be\label{eq:lambda}
\!\!\! \hat{\lambda}_\delta = \frac{\tilde{p}_2(0)}{p_1(0)} =\Phi(h)\!-\! \frac{(h\!+\!2\delta) \kappa_{h,\delta}
 \!+\!\varphi(h) [ \Phi(-3\delta)  e^{\delta^2/2-h^2/2-2\delta h }\! -\!  \Phi(h-\delta)  e^{-3\delta h-  7\delta^2/2}]
  }
 {(h+2\delta) \left[\Phi(h)-\Phi(-\delta) e^{-(h+\delta)(h+3\delta)/2}\right]} . \;\;
\ee
Moreover, $\tilde{p}_2(x)$ is a rather accurate approximation to the (non-normalized) eigenfunction $p(x)$ in \eqref{eq:lambdaD}. In Fig~\ref{d_and_l}(a) we have plotted $p_0(x)$, $p_1(x)$ and the uncorrected $p_1(x)$ obtained by letting $\rho \rightarrow 0$  for  particular $L$ and $h$.

An alternative way of approximating $\lambda_\delta$ and $p(x)$ from  \eqref{eq:lambdaD}
would be to use a methodology described in \cite{Quadrature} p.154 which is based on the Gauss-Legendre  discretization of  the interval $[-C,h]$, with some large $C>0$, into an $N$-point set $x_1, \ldots, x_N$ (the $x_i$'s are the roots of the \mbox{$N$-th} Legendre polynomial on $[-C,h]$), and the use of the Gauss-Legendre weights $w_i$ associated with points $x_i$; $\lambda_\delta$ and $p(x)$ are  then approximated by the largest eigenvalue and associated eigenvector of the matrix
$
D^{1/2}K_0D^{1/2},
$
where $D = \text{diag}({w}_i)$, and  $(K_0)_{i,j} = p_{h,\delta}(x_i\given x_j)$. If $N$ is large enough then the resulting approximation to $\lambda_\delta$ is arbitrarily accurate and we use it as the true $\lambda_\delta$ in our numerical comparisons. Numerical simulations show that the value of $c_1$ is not close enough to $\lambda_\delta$ but  $c_2$ is.
However, more interestingly, we see that $\hat{\lambda}_\delta$ defined by \eqref{eq:lambda} is very accurate and we suggest to use it because of its explicit form; this is demonstrated in Fig~\ref{d_and_l}(b), where $\hat{\lambda}_\delta$  (solid red line)
is visually indistinguishable from   ${\lambda}_\delta$ obtained using the Gauss-Legendre quadrature (dotted black line). In this figure, the dashed green line corresponds to the uncorrected ${\lambda}_\delta$ ($\rho \rightarrow 0$), where we once again see a very significant difference between the corrected and uncorrected approximations.

We have discussed how to correct $\lambda$ for discrete time. We shall now discuss item (a) of Section~\ref{Correct_t_gt_1}, which concerns correction of the BCP ${\cal P}_{\zeta}(1,h)$ for discrete time. For correcting  ${\cal P}_{\zeta}(1,h)$ we can routinely use $\delta$ as in \eqref{delta} but numerical results indicate that we get a better resulting approximation, especially for small $L$ and $M$,  if we use $\gamma = \rho_L /T^{1/4} $, so that the BCP ${\cal P}_{\mathbb{X}}(1,h)$ is approximated by ${{\cal P}}_{\zeta,\gamma}(1,{h})$. We believe that the fact that ${\cal P}_{\mathbb{X}}(1,h)$ $\cong$ ${{\cal P}}_{\zeta,\delta}(1,{h})$ is not accurate enough is due to the fact that the densities $p_i(x)$ are not exactly the densities of $\xi_i$.
To summarise, the CDA for $T>1$ is the following approximation.\\

\noindent{\bf Approximation 6.} \textit{For $M> L$ (that is, $T > 1$), the CDA for the BCP \eqref{eq:prob-xi} is}
\be\label{correct_diff}
{\cal P}_{\mathbb{X}}(M,h) \cong 1-\left[1-{{\cal P}}_{\zeta,\gamma}(1,{h})\right] \,\hat{\lambda}_\delta ^{{ T}-1},
\ee
\textit{where $\hat{\lambda}_\delta$ is given by \eqref{eq:lambda} and }
\begin{eqnarray*}
{{\cal P}}_{\zeta,\gamma}(1,{h}) = 1 - \Phi(h+\gamma)\,\Phi(h) + \frac{\varphi(h+\gamma)}{\gamma}\Phi(h) - \frac{\varphi(h)e^{-2h{\gamma}}}{\gamma}\Phi(h-\gamma)
\end{eqnarray*}
\textit{ with $\gamma = \rho_L /T^{1/4} = 0.5826/(L^{1/2}T^{1/4})$ }.

\subsection{Approximation by J. Glaz and coauthors}
\label{sec:Glaz}

The Glaz approximation for the BCP ${\cal P}_{\mathbb{X}}(M,h)$ (developed in \cite{Glaz_old,Glaz2012} and discussed in the Introduction) is as follows.\\

\noindent{\bf Approximation 7.} (Glaz approximation) For $M \ge 2L$ (so that $T=M/L \ge 2$)
\be
\label{eq:Glaz}
{\cal P}_{\mathbb{X}}(M,h) \cong 1- (1-{\cal P}_{\mathbb{X}}(2L,h) ) \left[
 \frac{1-{\cal P}_{\mathbb{X}}(2L,h) }{1-{\cal P}_{\mathbb{X}}(L,h) }
  \right]^{T-2} \, ,
\ee
where ${\cal P}_{\mathbb{X}}(2L,h) $ and ${\cal P}_{\mathbb{X}}(L,h) $ are evaluated using R algorithms for the multivariate normal distribution.

The approximation \eqref{eq:Glaz} is defined  for $M \ge 2L$ and  requires numerical  evaluation of $L+1$
and  $2L+1$   dimensional integrals (which are the BCP ${\cal P}_{\mathbb{X}}(L,h) $ and ${\cal P}_{\mathbb{X}}(2L,h)$ respectively)
using the so-called `GenzBretz' algorithm for numerical evaluation of multivariate normal probabilities, see \cite{genz2009computation,GenzR}.
  Whilst the accuracy of Approximation 7 is very high  and in fact very similar the accuracy of the CDA (Approximation 6), the nature of Approximation 7 results in high computational cost and run-time when compared to other approximations discussed in this paper (especially for large L);
  note also that different integrals should be computed for different values of $h$.
  Moreover, the `GenzBretz' algorithm uses Monte-Carlo simulations so that
  for reliable estimation of high-dimensional integrals (especially when $L$ is large) one needs to make a lot of averaging.

  We have not provided results of comparison of Approximation 7 with other approximations for the BCP as the accuracy of Approximations 6 and 7 was very close.
  Note also that there is  strong similarity between the forms of these two approximations.
  Indeed, from \eqref{correct_diff} we can write  the CDA in the form
\begin{eqnarray*}
1-(1-{{\cal P}}_{\zeta,\gamma}(1,{h})) \hat{\lambda}_\delta ^{{ T}-1}
= 1 - \underbrace{\hat{\lambda}_\delta(1-{{\cal P}}_{\zeta,\gamma}(1,{h}))}_{(a)}{\underbrace{{\hat{\lambda}_\delta}}_{(b)}}^{T-2},
\end{eqnarray*}
where the terms (a) and (b) are as  $(1-{\cal P}_{\mathbb{X}}(2L,h) )$ and
$({1-{\cal P}_{\mathbb{X}}(2L,h) })/({1-{\cal P}_{\mathbb{X}}(L,h) })
 $
in Approximation 7, respectively.

\subsection{Simulation study}\label{sec:sim-study2}

\begin{figure}[h]
\begin{center}
 \includegraphics[width=0.5\textwidth]{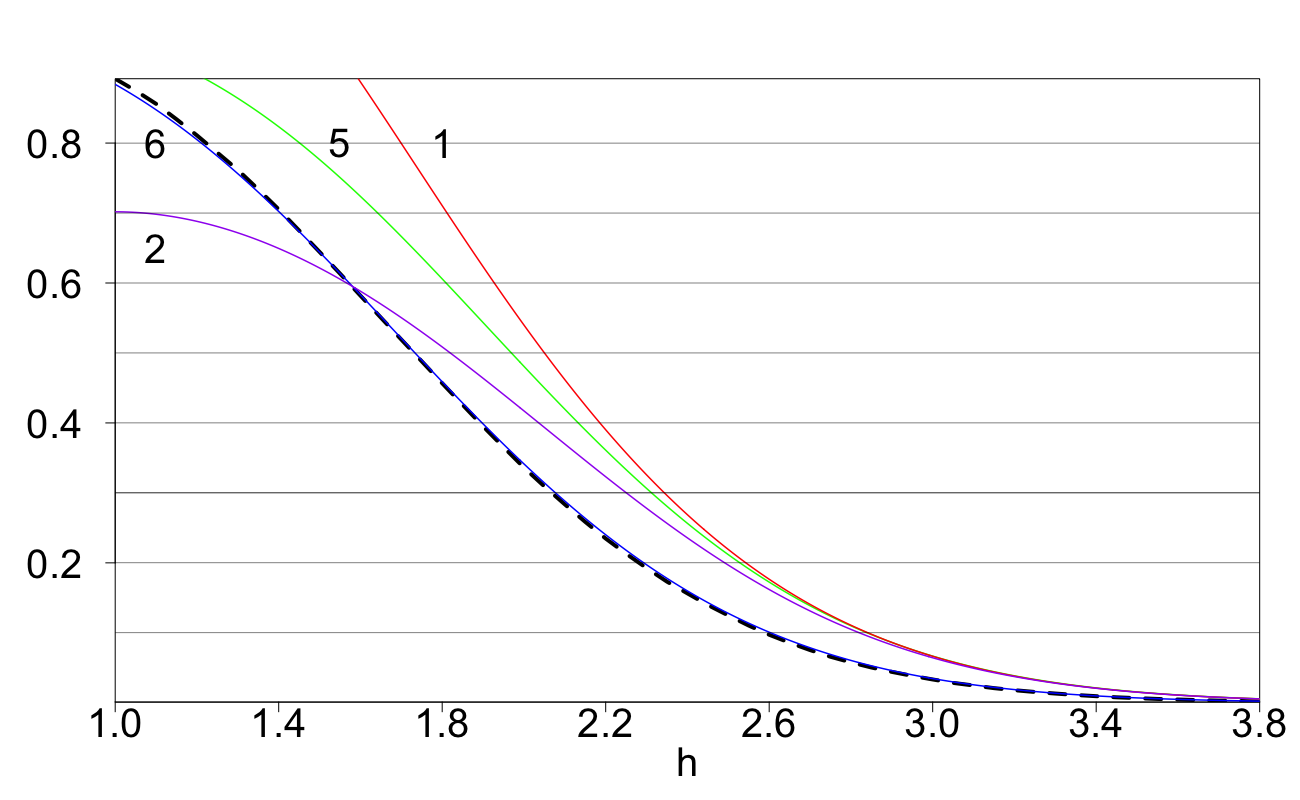}\includegraphics[width=0.5\textwidth]{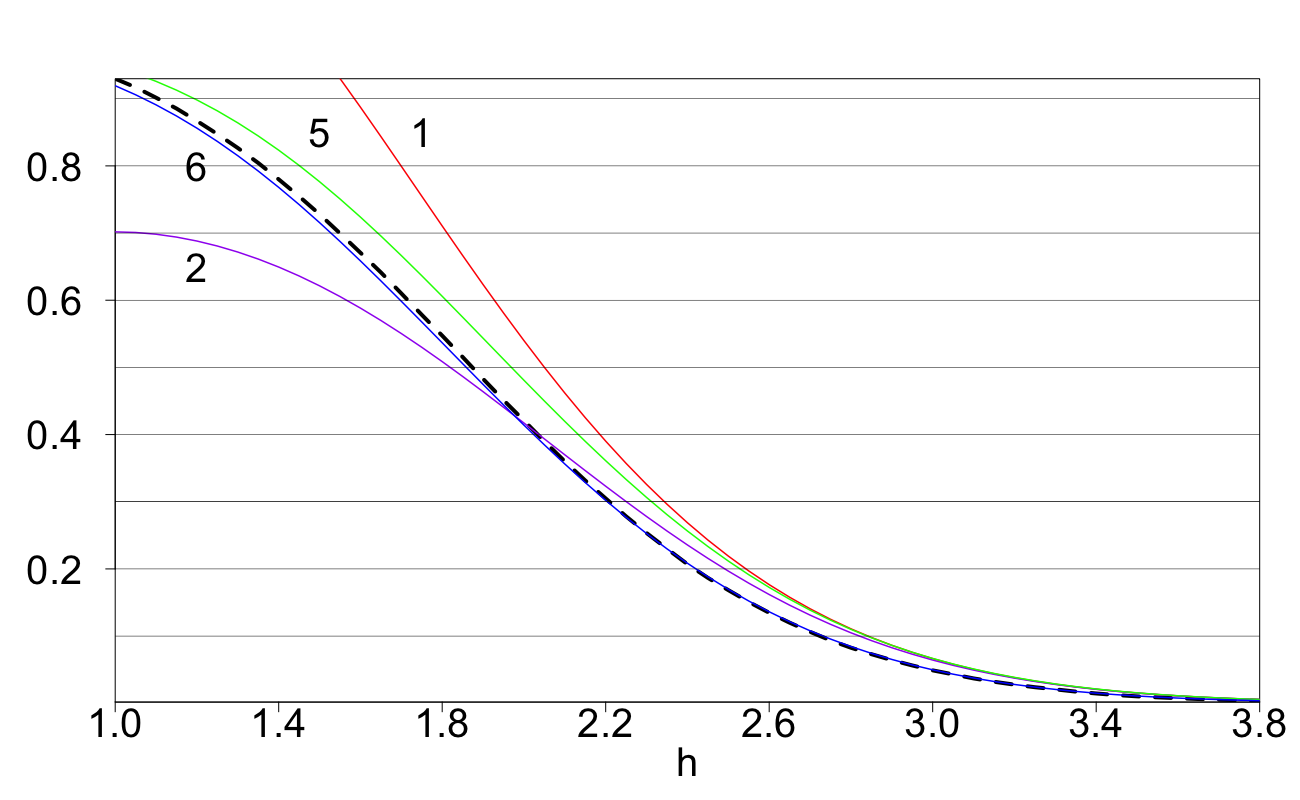}
\end{center}
\caption{Empirical probabilities of reaching the barrier $h$ and four approximations.
Left: $L=10$, $M=50$, $T=5$. Right: $L=50$, $M=250$, $T=5$. }
\label{L10_L50}
\end{figure}

\begin{figure}[h]
\begin{center}
 \includegraphics[width=0.5\textwidth]{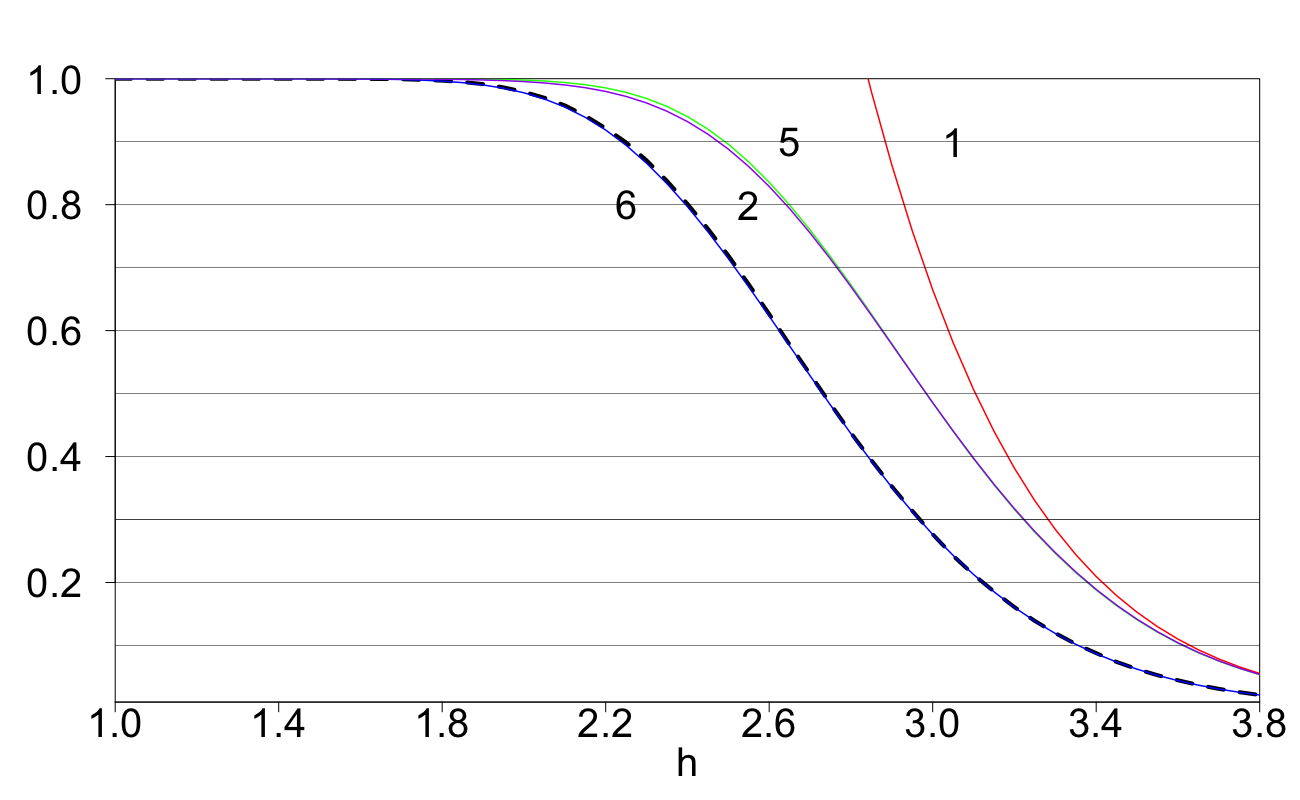}\includegraphics[width=0.5\textwidth]{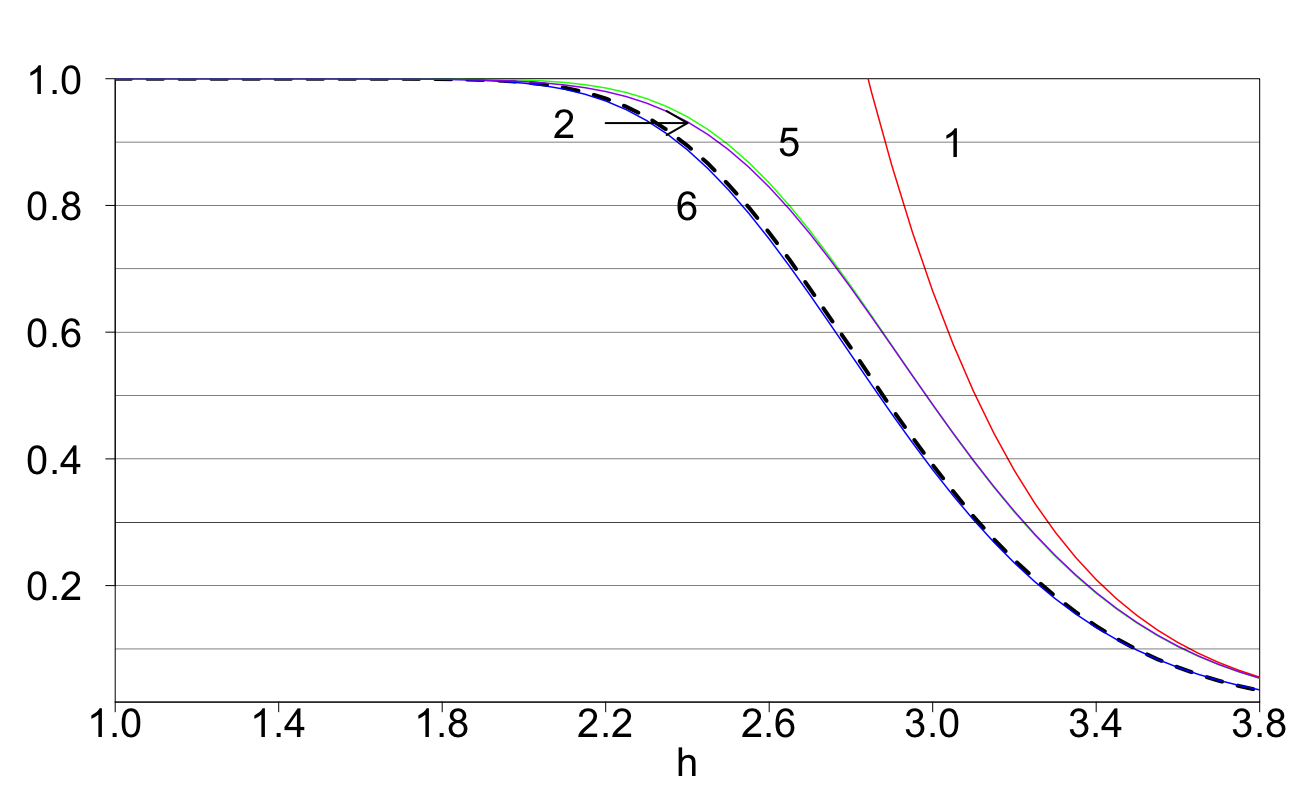}
\end{center}
\caption{Empirical probabilities of reaching the barrier $h$ and four approximations.
Left: $L=10$, $M=500$, $T=50$. Right: $L=50$, $M=2500$, $T=50$. }
\label{L10_L50_T_50}
\end{figure}

In this section we study the quality of the Durbin (Approximation 1), PCH (Approximation~2), Diffusion (Approximation 5) and CDA (Approximation 6) approximations for the BCP
${\cal P}_{\mathbb{X}}(M,h)$, defined in  \eqref{eq:prob-xi},  when $M> L$ (so that $T>1$).
The styles  of Fig.~\ref{L10_L50}, Fig.~\ref{L10_L50_T_50} and Table~\ref{Relative_error_table2} are exactly the same as of Fig.~\ref{Quad1}, Fig.~\ref{Quad2} and Table~\ref{Relative_error_table}, respectively,  and are described in
the beginning of Section~\ref{sec:sim-study}.
Similar to the case  $T \le 1$, we conclude that the CDA provides very accurate approximations and significantly
outperforms the Diffusion, Durbin and PCH approximations. Note also that for large $T$ the PCH approximation is very close to the Diffusion approximation;
this can be seen in Fig~\ref{L10_L50_T_50}, where (for $T=50$) these two approximations basically coincide for all~$h$.

\vspace{-0.5cm}
\begin{table}[h]
\centering
\caption{Relative error of CDA for given BCP}  \label{Relative_error_table2}     
\begin{tabular}{|c||c|c|c|c|c|c|}
\hline
BCP & $L\!=\!10,M\!\!=\!50$& $L\!=\!10,M\!\!=\!500$  & $L\!=\!50,M\!\!=\!250$& $L\!=\!50,M\!\!=\!2500$\\
\hline
0.05 & 0.596 \% & 0.028  \% &   0.133  \%&  0.054 \%  \\
0.10 &  0.657 \%  & 0.030  \% &  0.146   \%&  0.057  \% \\
0.15 &  0.455 \% & 0.031  \% &    0.390  \%& 0.208  \% \\
0.20 &  0.570 \% &  0.192 \% &  0.165  \%&  0.184  \%  \\
\hline
\end{tabular}
\end{table}

\section{Approximating ${\rm ARL}_h(\mathbb{X})$ }
\label{ARL_section}
As shown in Sections~\ref{sec:sim-study} and \ref{sec:sim-study2}, the CDA accurately approximates ${\cal P}_{\mathbb{X}}(M,h)$. The CDA has different forms depending on whether $M\le L$ or $M>L$,
see \eqref{Corrected_Diffusion} and \eqref{correct_diff} respectively. Thus, from \eqref{eq:first_pass_time}, the CDA
leads to the following approximation for  the probability density function of $\tau_{h}(\mathbb{X})$$/L$:
\begin{equation*}\label{q_t_h}
\hat{q}_h(t)=
\begin{cases}
\frac{d}{dt} \left \{ \int_{-\infty}^{h} Q_{h,\rho}(tL,x_0) \varphi(x_0) dx_0
\right\}, \,\,\,\, 0 < t \le 1\\[5pt]
\frac{d}{dt} \left\{ 1-\left[1-{{\cal P}}_{\zeta,\gamma}(1,{h})\right] \,\hat{\lambda}_\delta ^{{ t}-1}\,  \right\}, \,\,\,\, t >1.
\end{cases}
\end{equation*}
For $t >1$, one can easily get an explicit form of $\hat{q}_h(t)$. However, we were
 unable to obtain an explicit form of $\hat{q}_h(t)$ for $t < 1$  but this function can easily be numerically evaluated.
For  large ARL (and hence large $h$), the probability of exceeding $h$ in the interval $(0,1]$ is very small and the impact of  $\hat{q}_h(t)$ for $t < 1$ in the ARL approximation is minimal.

\begin{figure}[h]
\begin{center}
 \begin{subfigure}{0.45\linewidth} \centering
     \includegraphics[scale=0.17]{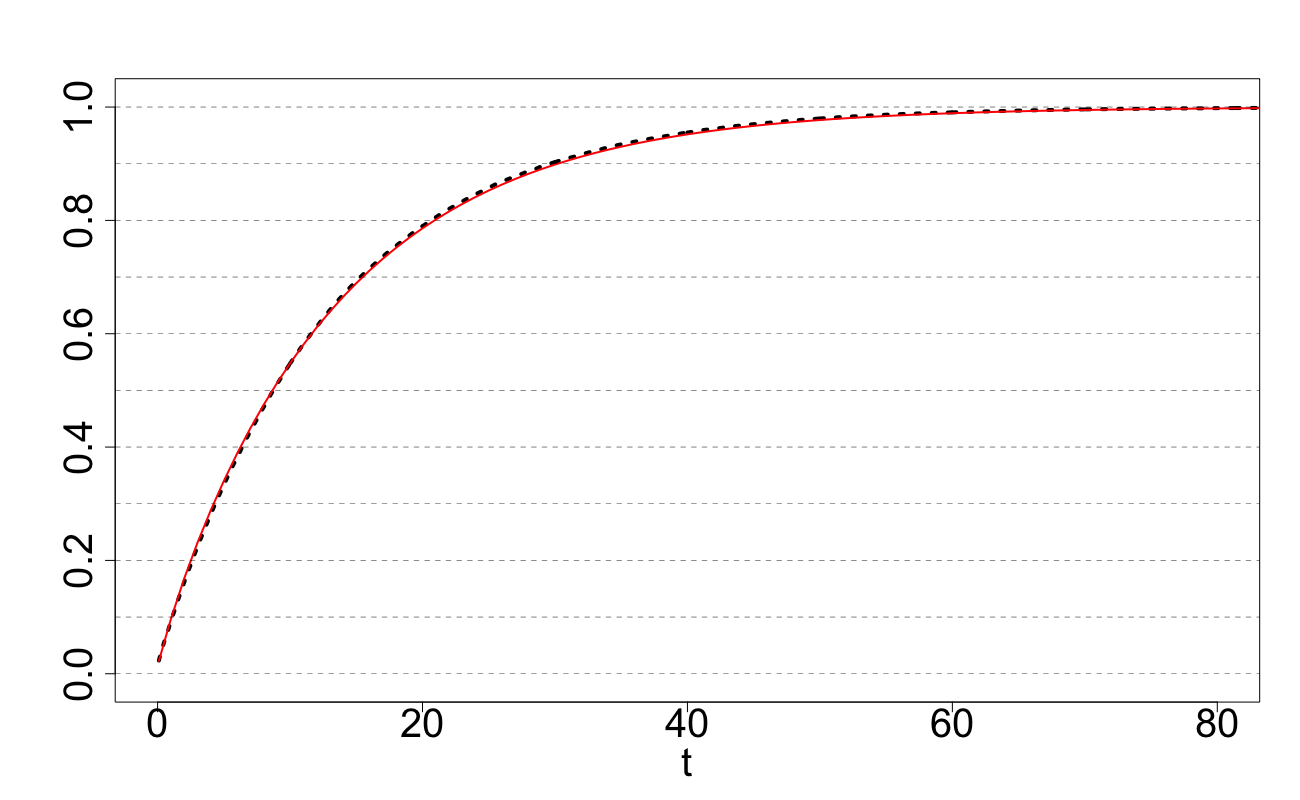}
     \caption{Left: $L=10$ and $h=2$ }
   \end{subfigure}
   \begin{subfigure}{0.45\linewidth} \centering
     \includegraphics[scale=0.17]{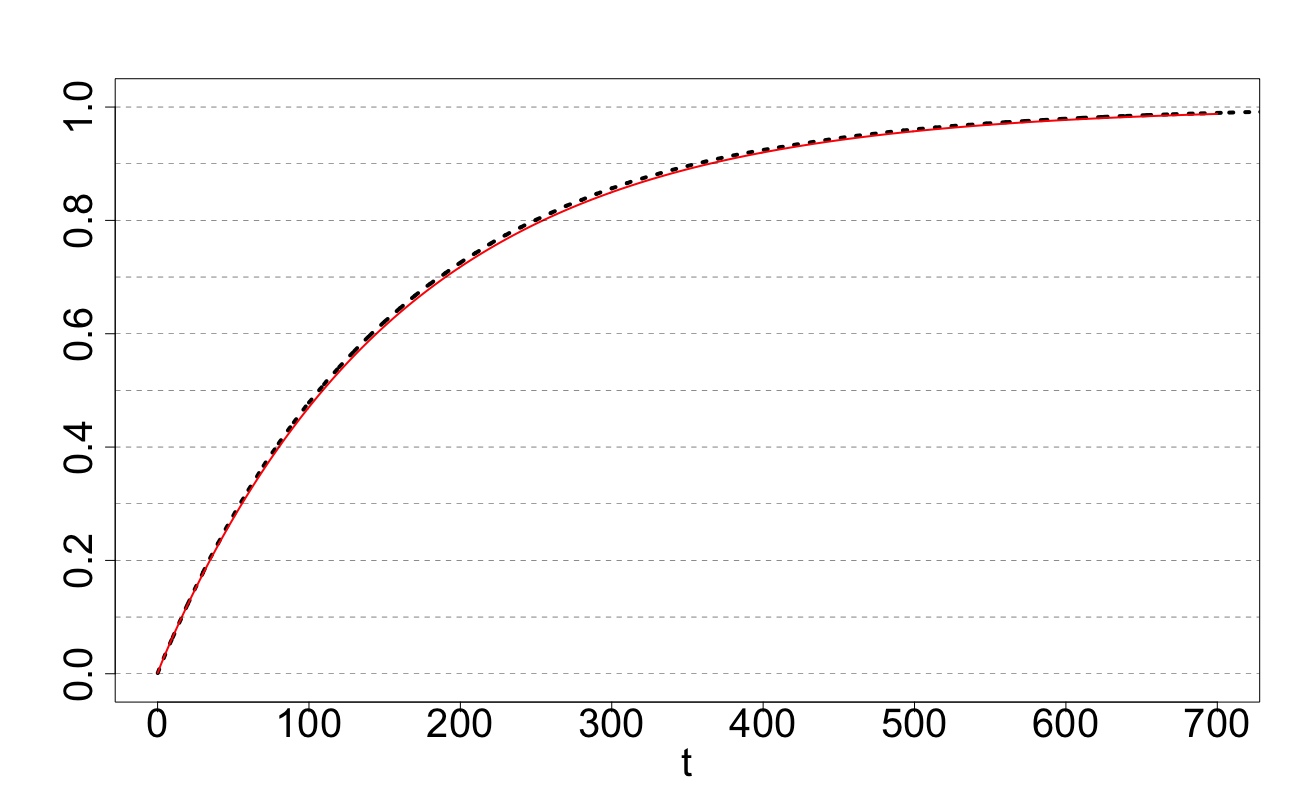}
     \caption{Right: $L=10$ and $h=3$}
   \end{subfigure}
\end{center}
\caption{ ${F}_h(t)$ and its approximation $\hat{F}_h(t)$ for $L=10$, $h=2$ and $3$.}
\label{Distribution1}
\end{figure}

\begin{figure}[h]
\begin{center}
 \begin{subfigure}{0.45\linewidth} \centering
     \includegraphics[scale=0.17]{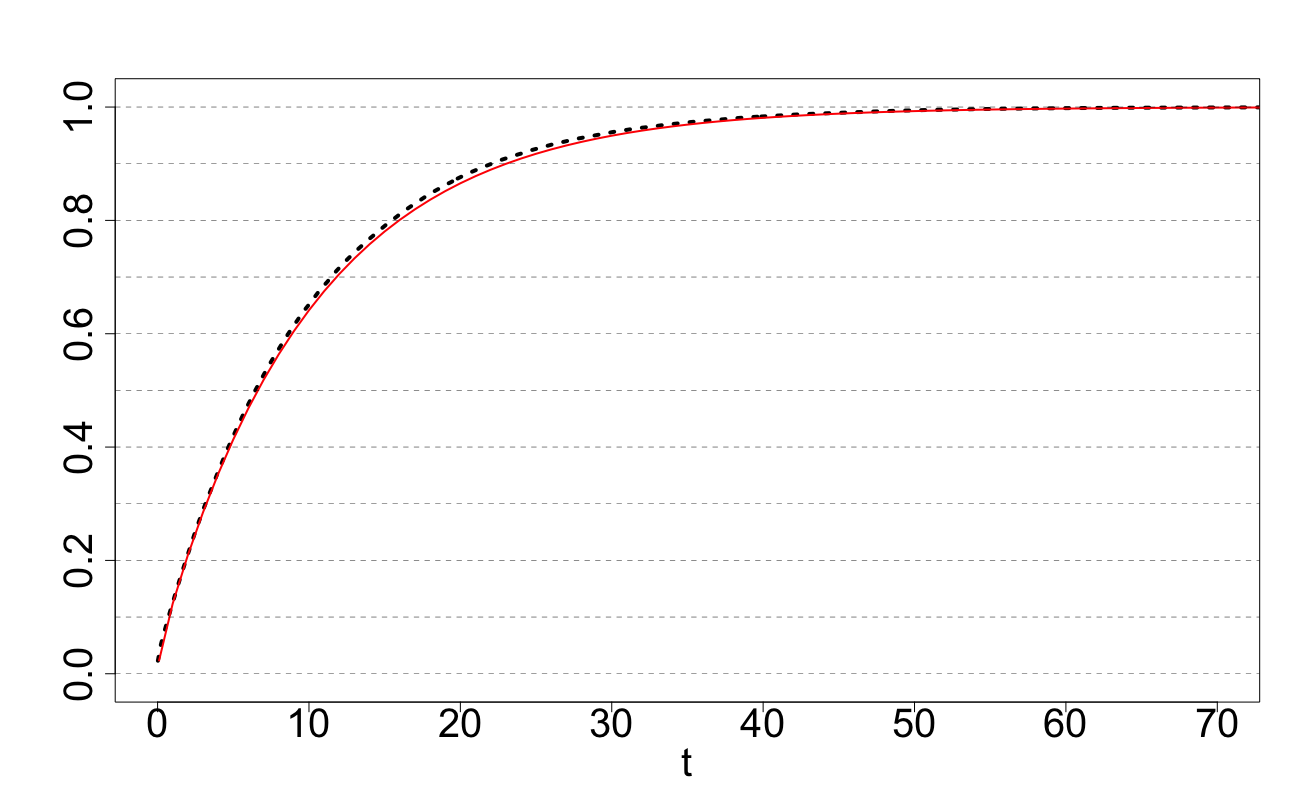}
     \caption{Left: $L=50$ and $h=2$ }
   \end{subfigure}
   \begin{subfigure}{0.45\linewidth} \centering
     \includegraphics[scale=0.17]{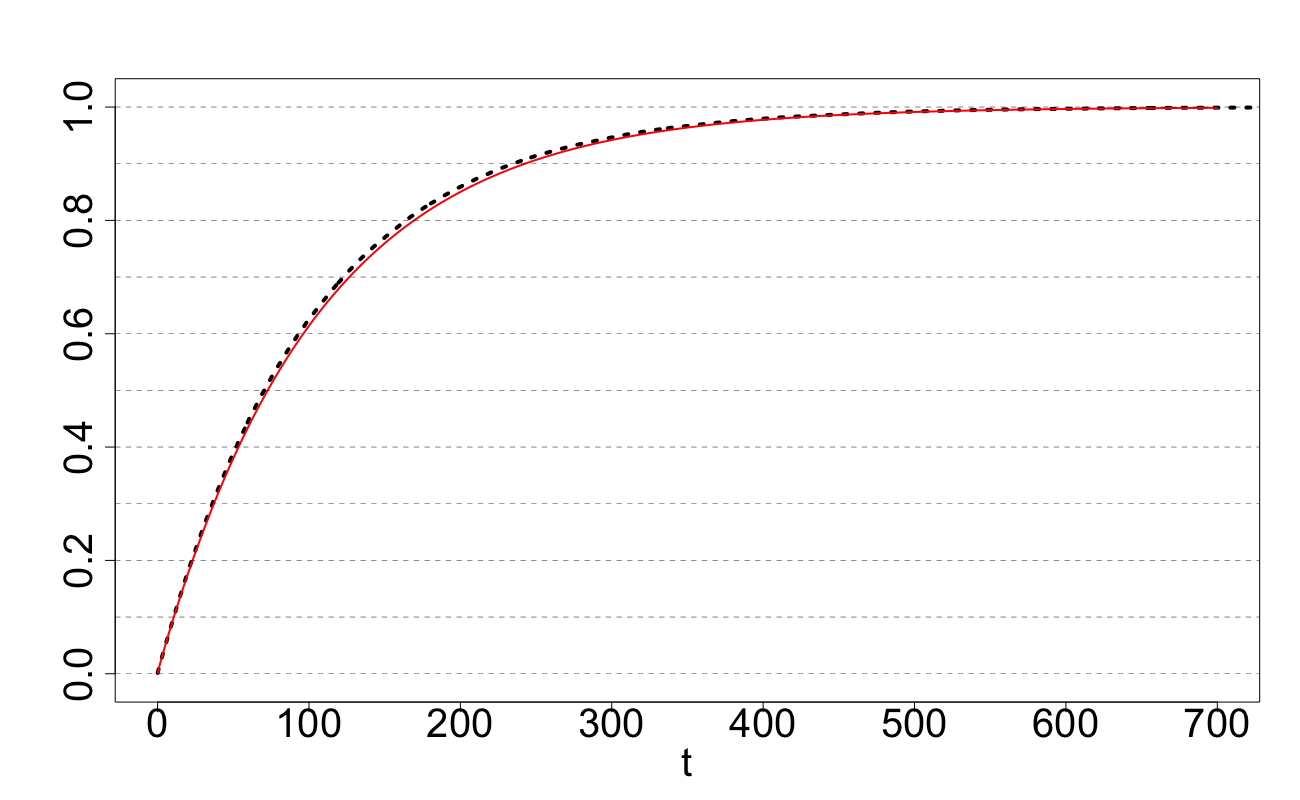}
     \caption{Right: $L=50$ and $h=3$}
   \end{subfigure}
\end{center}
\caption{${F}_h(t)$ and its approximation $\hat{F}_h(t)$ for $L=50$, $h=2$ and $3$. }
\label{Distribution2}
\end{figure}

Denote by ${F}_h(t)$ the true cumulative distribution function (c.d.f.) of $\tau_{h}(\mathbb{X})/L$.
The c.d.f. of the CDA of $\tau_{h}(\mathbb{X})/L$ is $\hat{F}_h(t)$ defined by
\be\label{CDF_approx}
\hat{F}_h(0)=1-\Phi(h), \;\; \hat{F}_h(t)= 1-\Phi(h) + \int_0^t \hat{q}_h(u) du\;\;\;{\rm for} \;t>0\,.
\ee
The accuracy of this approximation for a selection of parameter choices is demonstrated in Fig.~\ref{Distribution1} and Fig.~\ref{Distribution2}. The CDA for ${\rm ARL}_h(\mathbb{X})$  is
\begin{equation}
{\rm ARL}_h(\mathbb{X})  =  \E \tau_{h}(\mathbb{X}) \cong L \int_0^{\infty}{s \hat{q}(s,h) ds}
 \, .\label{ARL_formula}
\end{equation}

Fig.~\ref{Distribution1} and~\ref{Distribution2} demonstrate that   \eqref{CDF_approx} very  accurately  approximates the distribution of $\tau_{h}(\mathbb{X})$.

In this paper, we define ARL in terms of the number of random variables $\xi_{n}$ rather than number of random variables $\varepsilon_j$. This means we have to slightly modify the following Glaz approximation for ARL given in \cite{Glaz2012}, since such an approximation considers the number of random variables $\varepsilon_j$. This can be simply done by subtracting $L$ from the ARL approximation in  \cite{Glaz2012}.
From which, the Glaz approximation for ${\rm ARL}_h(\mathbb{X})$ is as follows:
\be \label{Glaz_ARL}
 \E _G(\tau_{h}(\mathbb{X})) = \sum_{j=L}^{2L}(1\!-\!{\cal P}_{\mathbb{X}}(j-L,h)) \!+\! \frac{1-{\cal P}_{\mathbb{X}}(L,h)}{{\cal P}_{\mathbb{X}}(2L,h)-{\cal P}_{\mathbb{X}}(L,h)}\sum_{j=1}^{L}(1-{\cal P}_{\mathbb{X}}(L+j,h))\, , \;\;\;
\ee
where $x=(1-{\cal P}_{\mathbb{X}}(2L,h) )/(1-{\cal P}_{\mathbb{X}}(L,h))$.

 In Table~\ref{expected_run_length}  we assess the accuracy of  the CDA approximation  \eqref{ARL_formula}
 and also Glaz approximation \eqref{Glaz_ARL}. In these tables, the values of ${\rm ARL}_h(\mathbb{X})$  have been calculated using $100,000$ simulations. Due to the Monte Carlo methods used to compute the Glaz approximation, we have presented the average of 20 iterations of \eqref{Glaz_ARL} as well as providing confidence intervals.


\begin{table}[!h]
\centering
\caption{Approximations for ${\rm ARL}_h(\mathbb{X})$: $L=10$ (top) and $L=50$ (bottom)}
\label{expected_run_length}
\begin{tabular}{|c||c|c|c|c|c|c|c|c|c|}
\hline
$h$ & 1 & 1.25 &1.5 &1.75 & 2 &2.25 & 2.5 &2.75 & 3  \\
\hline
 \eqref{ARL_formula}&21 & 32 &  49 &78 &  128 & 222&  403 & 774 & 1579        \\
\eqref{Glaz_ARL}  & 21  & 31 &48  & 77 & $126$ &218 $\pm$ 1 &  394 $\pm$ 2   &757 $\pm$ 5 &   1545 $\pm$ 18       \\
${\rm ARL}_h(\mathbb{X})$ & 22 &33 & 49 &78 & 127 &  218 & 397  & 758 & 1551       \\
\hline
\end{tabular}

\medskip

\centering

\begin{tabular}{|c||c|c|c|c|c|c|c|c|c|}
\hline
$h$ & 1 & 1.25 &1.5 &1.75 & 2 &2.25 & 2.5 &2.75 & 3  \\
\hline
 \eqref{ARL_formula}& 85 & 128 & 195  & 303 &  489  &  819 &  1440  & 2672  &     5256     \\
\eqref{Glaz_ARL}  & {82}   &  {123}  & {187}    &  {292} &  {474} \!$\pm$\!\! 1   &  {791}  $\pm$ 3&      {1393}  $\pm$ 8 &  {2597} $\pm$ 23  &   {5121}  $\pm$ 82        \\
${\rm ARL}_h(\mathbb{X})$ & 83 & 124 & 189 & 294 & 471 & 791   & 1392  &2590 &  5110     \\
\hline
\end{tabular}
\end{table}

     Table~\ref{expected_run_length}  shows that   for small $h$
      the approximations developed in this paper  are very accurate and are similar to the Glaz approximation.
      For a large $L=50$, \eqref{Glaz_ARL}  can be considered more accurate than \eqref{ARL_formula}. However using  \eqref{Glaz_ARL} for a large $L$ is computationally expensive and results in a long run-time, especially if results are averaged. Increasing $L$ has no impact on the computational cost and run time of \eqref{ARL_formula}.

\begin{acknowledgements}
The authors are grateful to our colleague Nikolai Leonenko for intelligent discussions and finding the reference \cite{Harr}, which is essential for the material of Section 4.2.
\end{acknowledgements}

\section*{Appendices}

\subsection*{Appendix A: Proof of Lemma \ref{lem:corr}}
As correlation is invariant under linear transformations,  ${\rm Corr}(S_{0,L},S_{k ,L}) = {\rm Corr}(\xi_0,\xi_k)$.
From the definition \eqref{eq:sumsq2} we have ${\rm Corr}(S_{0,L} , S_{k ,L})={\rm Corr}(S_{n,L},S_{n+k ,L})$. The sum $S_{k,L}$ can be represented as
$$S_{k,L}= S_{0,L} - \sum_{j=1}^{k} \varepsilon_j + \sum_{j=L+1}^{L+k} \varepsilon_j \, .
$$
Using this representation, we  obtain
$$
\text{Cov}(S_{0,L},S_{k,L}) = \underbrace{(\sigma^2L+\mu^2L^2)}_{{ \E S_{0,L}^2}} - k\sigma^2 - \underbrace{\mu^2L^2}_{ (\E S_{0,L})^2} = \sigma^2L-k\sigma^2\, .
$$
Dividing this by var$(S_{0,L})$, from \eqref{eq:mean_vars}, we obtain
$
{\rm Corr}(S_{0,L},S_{k,L})=1-{k}/{L}\,
$
in the case $k\le L$.  The case $k> L$ is obvious.

\subsection*{Appendix B: Derivation of Durbin approximation}

We shall initially show $R'(0+)=-1 \ne 0$. We have
$$
\left.\frac{\partial R_\zeta(t,s)}{\partial s}\right|_{s=t+}=R(0+).
$$
Using (\ref{eq:correlation}) and the fact that $\Delta ={1}/{L}$,
we have
$$
R'(0+)=\lim_{L\rightarrow\infty}\frac{R(\Delta)-R(0)}{\Delta}=
-\lim_{L\rightarrow\infty}\frac{L}{L}=-1.
$$

The Durbin approximation for $q(t,h,\zeta_t)$ can be written as
\bea
q(t,h,\zeta_t) \cong b_0(t,h)f(t,h)\, ,
\eea
where
$$
f(t,h)\!=\!\frac{1}{\sqrt{2\pi
R_{\zeta}(t,t)}}\,e^{-\frac{h^2(t)}{2R_{\zeta}(t,t)}}\,,\;\;
b_0(t,h)\!=\!\left.-\frac{h(t)}{R_{\zeta}(t,t)} \frac{\partial
R_{\zeta}(s,t)}{\partial s}\right|_{s=t+}\!-\frac{dh(t)}{dt}\, .
$$
In view of (\ref{eq:first_pass_time})
the related approximation for
the first passage probability ${\cal P}_{\zeta}(T,h) $ is
$$ {\cal P}_{\zeta}(T,h)  \cong  \int_{0}^{ T}{b_0(t,h)f(t,h)dt}\, .$$
In the case when the threshold $h(t)=h$ is constant, using  Lemma \ref{Durbin} we obtain
$$
b_0(t,h)=-hR'(0+)=h,\;\;\;\;
q(t,h,\zeta_t) \cong \frac{h}{\sqrt{2\pi}}\,e^{-h^2/2}
$$
and therefore we obtain the following approximation.
\be
{\cal P}_{\mathbb{X}}(M,h)  \cong  {\cal P}_{\zeta}(T,h)  \cong \frac{h{ T}}{\sqrt{2\pi}}e^{-h^2/2}. \nonumber
\ee

\subsection*{Appendix C: Derivation of \eqref{Corrected_Diffusion_explicit}}
As $M=L$,
$
 {\cal P}_{\zeta,\rho}(L,{h}) = \int_{-\infty}^{h}\left(1 - \Phi\left(b+\hat{a}\right) + e^{-2\hat{a}b}\Phi\left(b-\hat{a}\right) \right)\varphi(x_0) dx_0 + 1- \Phi(h).
$
Using the fact $\hat{a} = (h-x_0)/2 + \rho_L$ and  $b = (h+x_0)/2$, we obtain:
\bea
 {\cal P}_{\zeta,\rho}(M,{h})\!&=& \!1 - \!\!\int_{-\infty}^{h}\!\!\!\Phi\left(h+\rho_L \right)\varphi(x_0) dx_0 + \int_{-\infty}^{h}\!\! e^{-h^2/2 + x_0^2/2 -  \rho_Lh -\rho_Lx_0}\Phi\left(x_0\!-\!\rho_L \right)\varphi(x_0) dx_0  \\
  &=& 1 - \Phi(h+\rho_L )\Phi(h) + \varphi(h)e^{-\rho_L h}\int_{-\infty}^{h} e^{-\rho_L x_0} \int_{-\infty}^{x_0 - \rho_L} \varphi(z) dz \, dx_0.
 \eea
Making the substitution $k=z+\rho_L$ in the rightmost integral, we obtain
\bea
 \varphi(h)e^{-\rho_L h}\int_{-\infty}^{h} \int_{-\infty}^{x_0 }  e^{-\rho_L x_0} \varphi(k-\rho_L) dk \, dx_0.
\eea
By then changing the order of integration:
\bea
 \varphi(h)e^{-\rho_L h}\int_{-\infty}^{h} \int_{k}^{h}  e^{-\rho_L x_0} \varphi(k-\rho_L) dx_0 \, dk =  \frac{\varphi(h)e^{-\rho_L h}}{\rho_L}\int_{-\infty}^{h} (e^{-\rho_Lk} - e^{-\rho_Lh})\varphi(k-\rho_L) dk.
\eea
By expanding the brackets, we obtain:
\bea
 \frac{\varphi(h)e^{-\rho_L h}}{\rho_L}\int_{-\infty}^{h} (e^{-\rho_Lk} - e^{-\rho_Lh})\varphi(k-\rho_L) dk \!&=&\!\!
  \frac{\varphi(h)e^{-\rho_L h - \rho_L^2/2}}{\rho_L}\int_{\infty}^{h} \varphi(k) dk - \frac{\varphi(h)e^{-2\rho_L h }}{\rho_L}\int_{\infty}^{h} \varphi(k-\rho_L) dk \\
 &=& \frac{\varphi(h+\rho_L)}{\rho_L}\Phi(h) - \frac{\varphi(h)e^{-2\rho_L h }}{\rho_L} \Phi(h-\rho_L).
\eea
Thus we obtain the required:
\bea
 {\cal P}_{\zeta,\rho}(L,{h}) = 1 - \Phi(h+\rho_L )\Phi(h) +\frac{\varphi(h+\rho_L)}{\rho_L}\Phi(h) - \frac{\varphi(h)e^{-2h\rho_L }}{\rho_L} \Phi(h-\rho_L)\, .
\eea

\bibliographystyle{spmpsci}

\bibliography{changepoint}

\begin{thebibliography}{10}
\providecommand{\url}[1]{{#1}}
\providecommand{\urlprefix}{URL }
\expandafter\ifx\csname urlstyle\endcsname\relax
  \providecommand{\doi}[1]{DOI~\discretionary{}{}{}#1}\else
  \providecommand{\doi}{DOI~\discretionary{}{}{}\begingroup
  \urlstyle{rm}\Url}\fi

\bibitem{Aldous}
Aldous, D.: Probability {A}pproximations via the {P}oisson {C}lumping
  {H}euristic.
\newblock Springer Science \& Business Media (1989)

\bibitem{Chu}
Chu, C.S.J., Hornik, K., Kaun, C.M.: {MOSUM} tests for parameter constancy.
\newblock Biometrika \textbf{82}(3), 603--617 (1995)

\bibitem{Durbin}
Durbin, J.: The first-passage density of a continuous {G}aussian process to a
  general boundary.
\newblock Journal of Applied Probability pp. 99--122 (1985)

\bibitem{genz2009computation}
Genz, A., Bretz, F.: Computation of Multivariate Normal and t Probabilities.
\newblock Lecture Notes in Statistics. Springer-Verlag, Heidelberg (2009)

\bibitem{GenzR}
Genz, A., Bretz, F., Miwa, T., Mi, X., Leisch, F., Scheipl, F., Hothorn, T.:
  {mvtnorm}: Multivariate Normal and t Distributions (2018).
\newblock \urlprefix\url{https://CRAN.R-project.org/package=mvtnorm}.
\newblock R package version 1.0-8:
  \textit{`https://CRAN.R-project.org/package=mvtnorm'}

\bibitem{Glaz_old}
Glaz, J., Johnson, B.: Boundary crossing for moving sums.
\newblock Journal of {A}pplied Probability \textbf{25}(1), 81--88 (1988)

\bibitem{Glaz2012}
Glaz, J., Naus, J., Wang, X.: Approximations and inequalities for moving sums.
\newblock Methodology and Computing in Applied Probability \textbf{14}(3),
  597--616 (2012)

\bibitem{glaz2009scan2}
Glaz, J., Pozdnyakov, V., Wallenstein, S.: Scan Statistics: Methods and
  Applications.
\newblock Birkhäuser, Boston (2009)

\bibitem{Harr}
Harrison, J.: Brownian motion and stochastic flow systems.
\newblock John Wiley and Sons (1985)

\bibitem{Mehr}
Mehr, C., McFadden, J.: Certain properties of {G}aussian processes and their
  first-passage times.
\newblock Journal of the Royal Statistical Society. Series B (Methodological)
  \textbf{27}(3), 505--522 (1965)

\bibitem{Quadrature}
Mohamed, J., Delves, L.: Computational Methods for Integral Equations.
\newblock Cambridge University Press (1985)

\bibitem{MZ2003}
Moskvina, V., Zhigljavsky, A.: An algorithm based on {S}ingular {S}pectrum
  {A}nalysis for change-point detection.
\newblock Communications in Statistics---Simulation and Computation
  \textbf{32}(2), 319--352 (2003)

\bibitem{ReedSimon}
Reed, M., Simon, B.: Methods of Modern Mathematical Physics: Scattering theory
  Vol. 3.
\newblock Academic Press (1979)

\bibitem{Shepp66}
Shepp, L.: Radon-{N}ikodym derivatives of {G}aussian measures.
\newblock The Annals of Mathematical Statistics pp. 321--354 (1966)

\bibitem{Shepp71}
Shepp, L.: First passage time for a particular {G}aussian process.
\newblock The Annals of Mathematical Statistics \textbf{42}(3), 946--951 (1971)

\bibitem{Shepp76}
Shepp, L., Slepian, D.: First-passage time for a particular stationary periodic
  {G}aussian process.
\newblock Journal of Applied Probability pp. 27--38 (1976)

\bibitem{Sieg_book}
Siegmund, D.: Sequential Analysis: Tests and Confidence Intervals.
\newblock Springer Science \& Business Media (1985)

\bibitem{Sieg_paper}
Siegmund, D.: Boundary crossing probabilities and statistical applications.
\newblock The Annals of Statistics \textbf{14}(2), 361--404 (1986)

\bibitem{ZhK1988}
Zhigljavsky, A., Kraskovsky, A.: Detection of abrupt changes of random
  processes in radiotechnics problems.
\newblock St. Petersburg University Press (1988).
\newblock (in Russian)

\end{thebibliography}


\end{document}